\documentclass[12pt]{article}
\usepackage{graphics}
\usepackage{bbm}
\usepackage{times}
\usepackage{theorem}

\def\theoremname{Theorem}
\def\lemmaname{Lemma}
\def\definitionname{Definition}
\def\remarkname{Remark}
\def\examplename{Example}
\def\proofname{Proof}
\def\testcasename{Test Case}

\theoremstyle{plain}
\newtheorem{theorem}{\theoremname}[section]
\newtheorem{lemma}[theorem]{\lemmaname}
\newtheorem{definition}[theorem]{\definitionname}
\newtheorem{remark}[theorem]{\remarkname}

\theorembodyfont{\rm}
\theoremheaderfont{\it}
\newtheorem{testcase}[theorem]{\testcasename}
\newtheorem{proof}{\proofname}
\newtheorem{numproof}[theorem]{\proofname}

\newcommand{\intlabel}[1]{\label{In \thetheorem: #1}}
\newcommand{\intref}[1]{\ref{In \thetheorem: #1}}

\newcommand{\setN}{\mathbbm{N}}
\newcommand{\setQ}{\mathbbm{Q}}
\newcommand{\setR}{\mathbbm{R}}
\newcommand{\setZ}{\mathbbm{Z}}
\newcommand{\ind}{1}

\newcommand{\Lip}{\mathop{\mathrm{Lip}}}

\renewcommand{\d}[1]{\,\mathrm{d#1}}
\newcommand{\pprime}{{\prime\prime}}

\newcommand{\sgn}{\mathrm{sign}}

\newcommand{\conv}{\mathrm{conv}}

\newcommand{\esslim}{\mathop{\mathrm{ess\,lim}}}
\newcommand{\esslimsup}{\mathop{\mathrm{ess\,lim\,sup}}}

\def\next{\\\noalign{\vspace{\jot}}}

\begin{document}

\title{Well-Balanced Schemes for the Initial Boundary Value Problem for 1D Scalar Conservation Laws}
\author{M.~Nolte \and D.~Kr\"oner}
\date{}
\maketitle

\begin{abstract}
We consider well-balanced numerical schemes for the following 1D scalar conservation law with source term:
\begin{displaymath}
\partial_t u + \partial_x f\left( u \right) + z^\prime\left( x \right) \, b\left( u \right) = 0.
\end{displaymath}
More precisely, we are interested in the numerical approximation of the initial boundary value problem for this equation. While our main concern is a convergence result, we also have to extend Otto's notion of entropy solutions to conservation laws with a source term. To obtain uniqueness, we show that a generalization, the so-called entropy process solution (see \cite{eymard}), is unique and coincides with the entropy solution.

If the initial and boundary data are in $L^\infty$, we can establish convergence to the entropy solution. Showing that the numerical solutions are bounded we can extract a weak$\ast$-convergent subsequence. Identifying its limit as an entropy process solution requires some effort as we cannot use Kru\v zkov-type entropy pairs here. We restrict ourselves to the Engquist-Osher flux and identify the numerical entropy flux for an arbitrary entropy pair. By the uniqueness result, the scheme then approximates the entropy solution and a result by Vovelle then guarantees that the convergence is strong in $L^p$, $1 \le p < \infty$.
\end{abstract}
\vfill

\section{Introduction}

In this paper we will consider well-balanced numerical schemes for the initial boundary value problem
\begin{equation} \label{ibvp}
\begin{array}{rcll}
\displaystyle \partial_t u + \partial_x f\left( u \right) + z^\prime\left( x \right) b\left( u \right) &=& 0 & \enspace \mbox{in} \enspace \Omega_T,\next
\displaystyle u\left( \cdot, 0 \right) &=& u_0 &\enspace \mbox{in} \enspace \Omega,\next
\displaystyle u\left( x_l, \cdot \right) &=& u_l &\enspace \mbox{on} \enspace \left] 0, T \right[, \enspace \mbox{if} \enspace f^\prime\left( u \right) > 0,\next
\displaystyle u\left( x_r, \cdot \right) &=& u_r &\enspace \mbox{on} \enspace \left] 0, T \right[, \enspace \mbox{if} \enspace f^\prime\left( u \right) < 0,
\end{array}
\end{equation}
where $\Omega := \left] x_l, x_r \right[$, $\Omega_T := \Omega \times \left] 0, T \right[ $ and
\begin{equation} \label{ass}
\begin{array}{c}
\displaystyle
u_0 \in  L^\infty\left( \Omega \right),\quad
u_r, u_l \in L^\infty\left( \left] 0, T \right[ \right),\quad
f \in C^2( \setR, \setR ),\next
\displaystyle
z\in H^{1,\infty}\left( \Omega \right),\quad
b\in C^1\left( \setR \right),\enspace
b^\prime \in L^\infty\left( \setR \right)
\end{array}
\end{equation}
are given data. Furthermore we define
\begin{equation} \label{definition D}
D\left( s \right) := \int_0^s \frac{f^\prime\left( \xi \right)}{b\left( \xi \right)} \d\xi
\end{equation}
and assume that
\begin{equation} \label{assumptions-D}
D \in C^1\left( \setR \right),\quad
D\left( \setR \right) = \setR,\quad
\inf_\setR D^\prime > 0.
\end{equation}

For the corresponding initial value problem Greenberg et al.\ \cite{greenberg} have developed and investigated a well-balanced scheme which is much more efficient than standard schemes. In particular stationary solutions can be approximated by these schemes with less numerical investment than by the classical schemes. Further results including error estimates and numerical experiments for the initial value problem have been obtained by Gosse \cite{gosse}, Perthame \cite{perthame}. Corresponding results in multiple space dimensions can be found in \cite{botchorishvili} and \cite{pironneau}.

If $b=0$ there are many results for the initial boundary value problem concerning existence of entropy solutions (see \cite{bardos}, \cite{otto}) and convergence of numerical solutions to the entropy solution \cite{vovelle}. Now in this paper we want to analyse the combination of both, i.e.\ well-balanced schemes for the case $b \neq 0$.

The outline of this paper is as follows. In Section \ref{section-process-sol} we will repeat the definition of entropy and entropy process solutions. The numerical scheme will be described in Section \ref{well}. The main result and its proof will be given in Section \ref{main}. Finally in Section \ref{numerics} we will show the advantages of these well-balanced schemes in some numerical experiments.

\section{Entropy Process Solutions} \label{section-process-sol}

In \cite{otto} Otto has defined boundary entropy pairs for problems without source terms, which we will generalize here to problems with source terms. As in \cite{otto} we will define boundary entropy pairs as follows.

\begin{definition}[Boundary Entropy Pairs] \label{entropy}
\begin{enumerate}
\item Let $\eta \in C^2( \setR )$ be convex, $q \in C^1( \setR )$ and $q^\prime = \eta^\prime f^\prime$. Then $( \eta, q )$ is called an entropy pair for the partial differential equation in (\ref{ibvp}).

\item Let $H \in C^2( \setR^2 )$, $Q \in C^1( \setR^2 )$, $( H( \cdot, w ), Q( \cdot, w ) )$ be an entropy pair for all $w \in \setR$ so that
\begin{displaymath}
H\left( w, w \right) = \partial_1 H\left( w, w \right) = Q\left( w, w \right) = 0 \quad \forall w \in \setR.
\end{displaymath}
Then $( H, Q )$ is called a boundary entropy pair in the sense of Otto for the partial differential equation in (\ref{ibvp}).

\item Let $( \eta, q )$ be an entropy in the sense of (1) and $w\in \setR$ such that
\begin{displaymath}
\eta\left( w \right) = \eta^\prime\left( w \right) = q\left( w \right) = 0.
\end{displaymath}
Then $( \eta, q )$ is called a boundary entropy pair for the partial differential equation  in (\ref{ibvp}).
\end{enumerate}
\end{definition}

There is a close and simple relation between the boundary entropy pairs in the sense of (2) and (3) in \definitionname~\ref{entropy}.
\begin{remark}
\begin{enumerate}
\item If $( H, Q )$ is a boundary  entropy pair in the sense of Otto (see \definitionname~\ref{entropy}(2)) and $w \in \setR$ then
\begin{displaymath}
\eta\left( u \right) := H\left( u, w \right), \quad q\left( u \right) := Q\left( u, w \right)
\end{displaymath}
form a boundary entropy pair (see \definitionname~\ref{entropy}(3)).

\item If $( \eta, q )$ is a boundary entropy pair (see \definitionname~\ref{entropy}(3)) and $w\in \setR$ such that
$\eta( w ) = \eta^\prime( w ) = q( w ) = 0$ then
\begin{displaymath}
H\left( u, v \right) := \eta\left( u - \left(v - w\right) \right),
\quad Q\left( u, v \right) := q\left( u - \left(v - w\right) \right)
\end{displaymath}
form a boundary  entropy pair in the sense of Otto (see \definitionname~\ref{entropy}(2)).
\end{enumerate}
\end{remark}

Now we use the boundary entropy pairs from \definitionname~\ref{entropy}(3) to define an entropy solution of the initial boundary value problem (\ref{ibvp}).
\begin{definition}[Entropy Solution] \label{def-entropy-sol}
Let $u\in L^\infty( \Omega_T )$, $C := \Vert u \Vert_{L^\infty( \Omega_T )}$ and for all boundary entropy pairs $( \eta, q )$ and all $\varphi \in C_0^\infty( \overline\Omega \times [ 0, T [ )$, $\varphi \ge 0$, let
\begin{equation} \label{entropy-sol}
\begin{array}{c}
\displaystyle
\int_0^T \int_\Omega \eta\left( u \right) \partial_t \varphi + q\left( u \right) \partial_x \varphi
    - \eta^\prime\left( u \right) z^\prime\left( x \right) b\left( u \right) \varphi \d x \d t\next
\displaystyle
+ \int_\Omega \eta\left( u_0( x ) \right) \varphi( x, 0 ) \d x\next
\displaystyle
+ \Lip_{\left[ -C, C \right]}\left( f \right) \int_0^T \eta\left( u_r\left( t \right) \right) \varphi\left( x_r, t \right)
    + \eta\left( u_l\left( t \right) \right) \varphi\left( x_l, t \right) \d t
\ge 0,
\end{array}
\end{equation}
where
\begin{displaymath}
\Lip_{\left[ -C, C \right]}\left( f \right) := \sup_{-C \le u, v \le C}\frac{\left\vert f\left( u \right) - f\left( v \right)\right\vert}{\left\vert u - v \right\vert}.
\end{displaymath}
Then $u$ is called an entropy solution of (\ref{ibvp}). If $z^\prime = 0$ this definition corresponds to the definition of entropy solution by Otto (see \cite{malek}).
\end{definition}

The basic idea of the convergence proof for the numerical solutions consists in replacing the exact solution $u$ in (\ref{entropy-sol}) by the numerical solution $u_{\Delta x}$ as defined in \definitionname~\ref{def-scheme}. First we will prove that $u_{\Delta x}$ remains bounded for ${\Delta x} \to 0$ and then we have to control the limit in (\ref{entropy-sol}) for ${\Delta x}\to 0$. This can be easily done on the basis of the following lemma.

\begin{lemma}[Eymard, Gallou\"et, Herbin \cite{eymard}] \label{lemma-eymard}
Let $( u_k )_{k\in \setN} \subset L^\infty( \Omega_T )$ be a bounded sequence. Then there exists a subsequence $( u_{k^\prime} )_{k^\prime\in \setN}$ and a function $\mu \in L^\infty( \Omega_T \times ]0,1[ )$ such that for any $g \in C^0( \setR )$ we have
\begin{displaymath}
g\left( u_{k^\prime} \right) \to \int_0^1 g\left( \mu\left( x, t, \alpha \right) \right) \d\alpha \quad \mbox{weak}\ast.
\end{displaymath}
\end{lemma}

This gives us the motivation for the following definition of the entropy process solution.

\begin{definition}[Entropy Process Solution, see \cite{eymard}] \label{def-process-sol}
Let $\mu = \mu( x, t, \alpha )$ be an $L^\infty( \Omega_T \times ]0,1[ )$-function, $C := \Vert \mu \Vert_{L^\infty( \Omega_T \times ] 0, 1 [ )}$ and for all boundary entropy pairs $( \eta, q )$ and all $\varphi \in C_0^\infty( \overline\Omega \times [ 0, T [ )$, $\varphi \ge 0$, let
\begin{equation} \label{process-sol}
\begin{array}{c}
\displaystyle
\int_0^T \int_\Omega \int_0^1 \eta\left( \mu \right) \partial_t \varphi + q\left( \mu \right) \partial_x \varphi
    - \eta^\prime\left( \mu \right) z^\prime\left( x \right) b\left( \mu \right) \varphi \d\alpha \d x \d t\next
\displaystyle
+ \int_\Omega \eta\left( u_0( x ) \right) \varphi( x, 0 ) \d x\next
\displaystyle
+ \Lip_{\left[ -C, C \right]}\left( f \right) \int_0^T \eta\left( u_r\left( t \right) \right) \varphi\left( x_r, t \right)
    + \eta\left( u_l\left( t \right) \right) \varphi\left( x_l, t \right) \d t
\ge 0.
\end{array}
\end{equation}
Then $\mu$ is called an entropy process solution of (\ref{ibvp}).
\end{definition}

Now we are going to prove that the entropy process solution is unique. Later in Section \ref{main} we will show that the numerical solutions converge to an entropy process solution. Then the uniqueness will give us the existence of an entropy solution.

\begin{theorem}[Uniqueness of the Entropy Process Solution] \label{uniqueness}
Let $\mu, \nu \in L^\infty( \Omega_T \times ] 0, 1 [ )$ be two entropy process solutions of (\ref{ibvp}) with respect to the inital data $u_0 \in L^\infty( \Omega )$ and to the boundary data $u_l, u_r \in L^\infty( ] 0, T [ )$. Then $\mu = \nu$ a.e. on $\Omega_T \times ] 0, 1 [$ and
\begin{displaymath}
u\left( x, t \right) := \int_0^1 \mu\left( x, t, \alpha \right) \d\alpha
\end{displaymath}
is an entropy solution of (\ref{ibvp}).
\end{theorem}

Notice that this theorem also proves uniqueness of the entropy solution.

In order to prepare the proof of this theorem we need the following definition and some lemmata.

\begin{definition}[Semi-Kru\v zkov Entropy Pairs] \label{semi}
For $k\in \setR$ define the semi-Kru\v zkov entropy pairs
\begin{displaymath}
\begin{array}{rclrcl}
\displaystyle \eta_k^+\left( s \right) &:=& \displaystyle\left( s-k \right)^+,
    & \displaystyle q_k^+\left( s \right) &:=& \displaystyle \sgn^+\left( s-k \right) \left( f\left( s \right) - f\left( k \right) \right)\next
\displaystyle \eta_k^-\left( s \right) &:=& \displaystyle \left( s-k \right)^-,
    & \displaystyle q_k^-\left( s \right) &:=& \displaystyle \sgn^-\left( s-k \right) \left( f\left( s \right) - f\left( k \right) \right)
\end{array}
\end{displaymath}
where $s^+ := \sgn^+( s ) s$, $s^- := \sgn^-( s ) s$ and
\begin{displaymath}
\sgn^+\left( s \right) =
\left\lbrace\begin{array}{ll}
1 & \enspace \mbox{for} \enspace s > 0,\\
0 & \enspace \mbox{for} \enspace s \le 0,
\end{array}\right.
\quad \sgn^-\left( s \right) =
\left\lbrace\begin{array}{ll}
0 & \enspace \mbox{for} \enspace s \ge 0,\\
-1 & \enspace \mbox{for} \enspace s < 0.
\end{array}\right.
\end{displaymath}
\end{definition}

It can easily be shown that for any semi-Kru\v zkov entropy pair $( \eta, q )$ and any $u \in [ -C, C ]$ we have
\begin{displaymath}
\left\vert q\left( u \right) \right| \le \Lip_{\left[ -C, C \right]}\left( f \right) \eta\left( u \right).
\end{displaymath}

Now it turns out that we can use the semi-Kru\v zkov entropy pairs instead of the boundary entropy pairs in (\ref{process-sol}). This will be made precise in the next lemma.
\begin{lemma} \label{lemma-proc}
Let $\mu \in L^{\infty}( \Omega_T \times ] 0, 1 [ )$, $\mu = \mu( x, t, \alpha )$ be an entropy process solution of (\ref{ibvp}), $C := \Vert \mu \Vert_{L^\infty( \Omega_T \times ] 0, 1 [ )}$, $k\in \setR$ and $( \eta_k^+, q_k^+ )$, $( \eta_k^-, q_k^- )$ be semi-Kru\v zkov entropy pairs in the sense of \definitionname~\ref{semi}. Then we can replace $( \eta, q )$ in (\ref{process-sol}) by $( \eta_k^+, q_k^+ )$ and $( \eta_k^-, q_k^- )$ respectively, i.e.\ for all $\varphi \in C_0^\infty( \overline\Omega \times [ 0, T [ )$, $\varphi \ge 0$, we have
\begin{equation} \label{process-sol-kruzkov}
\begin{array}{c}
\displaystyle
\int_0^T \int_\Omega \int_0^1 \eta_k^\pm\left( \mu \right) \partial_t \varphi
    + q_k^\pm\left( \mu \right) \partial_x \varphi \d\alpha \d x \d t\next
\displaystyle
- \int_0^T \int_\Omega \int_0^1 \sgn^\pm\left( \mu - k \right) z^\prime\left( x \right) b\left( \mu \right) \varphi \d\alpha \d x \d t\next
\displaystyle
+ \int_\Omega \eta_k^\pm\left( u_0\left( x \right) \right) \varphi\left( x, 0 \right) \d x\next
\displaystyle
+ \Lip_{\left[ -C, C \right]}\left( f \right) \int_0^T \eta_k^\pm\left( u_r\left( t \right) \right) \varphi\left( x_r, t \right)
    + \eta_k^\pm\left( u_l\left( t \right) \right) \varphi\left( x_l, t \right) \d t
\ge 0.
\end{array}
\end{equation}
\end{lemma}

\begin{proof}
Let $p_{k,\delta}\left( x \right) := -\frac{1}{2 \delta^3}\left( x-k \right)^4 + \frac{1}{\delta^2}\left( x-k \right)^3$ and define
\begin{eqnarray*}
\eta_{k,\delta}^+\left( x \right) &:=& \left\lbrace
\begin{array}{ll}
0 & x \le k\\
p_{k,\delta}\left( x \right) & k \le x \le k+\delta\\
\left( x - k \right) - \frac{1}{2} \delta & x \ge k+\delta
\end{array}
\right.,\\
q^+_{k,\delta}\left( x \right) &:=& \int_k^x {\eta^+_{k,\delta}}^\prime\left( \xi \right) f^\prime\left( \xi \right) \d \xi.
\end{eqnarray*}
Then $( \eta^+_{k,\delta}, q^+_{k,\delta} )$ is a boundary entropy pair and we have uniformly in x
\begin{eqnarray*}
\bigl\vert \eta^+_{k,\delta}\left( x \right) - \eta^+_k\left( x \right) \bigr\vert &\le& \delta,\\
\bigl\vert {\eta^+_{k,\delta}}^\prime\left( x \right) - \sgn^+\left( x - k \right) \bigr\vert &\to& 0,\\
\bigl\vert q^+_{k,\delta}\left( x \right) - q_k\left( x \right) \bigr\vert &\le& \delta \left\Vert f^\prime \right\Vert_{L^\infty\left( [k,k+\delta] \right)}.
\end{eqnarray*}
Since $( \eta^+_{k,\delta}, q^+_{k,\delta} )$ is a boundary entropy pair we can put it into (\ref{process-sol}). For $\delta \to 0$ we obtain (\ref{process-sol-kruzkov}).
In a similar way we can prove (\ref{process-sol-kruzkov}) for $( \eta_k^-, q_k^- )$.
\end{proof}

The following lemma concerns the trace of an entropy process solution and follows the ideas of Otto (see \cite{otto}).

\begin{lemma} \label{trace}
Let $\mu \in L^\infty( \Omega_T \times ] 0, 1 [ )$, $\mu = \mu( x, t, \alpha )$, be an entropy process solution of (\ref{ibvp}). Then for all $v \in L^\infty( ] 0, T [ )$ and for all $\beta \in L^1( ] 0, T [ )$, $\beta \ge 0$ a.e., we have:
\begin{displaymath}
\begin{array}{l}
\displaystyle \esslim_{x \uparrow x_r} \int_0^T \int_0^1 \sgn^\pm\left( \mu - v \right) \left( f\left( \mu \right) - f\left( v \right) \right) \beta \d \alpha \d t\next
\ge \displaystyle - \Lip_{\left[ -C, C \right]}\left( f \right) \int_0^T \left( u_r\left( t \right) - v \right)^\pm \beta \d t,\next
\displaystyle \esslim_{x \downarrow x_l} \int_0^T \int_0^1 \sgn^\pm\left( \mu - v \right) \left( f\left( \mu \right) - f\left( v \right) \right) \beta \d \alpha \d t\next
\le \displaystyle \Lip_{\left[ -C, C \right]}\left( f \right) \int_0^T \left( u_l\left( t \right) - v \right)^\pm \beta \d t.
\end{array}
\end{displaymath}
The essential limits exist.
\end{lemma}

\begin{proof}
Let $w \in \setQ$ be fixed, $( \eta_w^+, q_w^+ ) $ be a semi-Kru\v zkov entropy pair and let
$\beta \in C_0^\infty( ] 0, T [ )$, $\beta \ge 0$, and $\psi \in C_0^\infty( \Omega )$, $\psi \ge 0$.
Then we obtain with \lemmaname~\ref{lemma-proc}:
\begin{equation} \intlabel{pre-monotony}
\begin{array}{rcl}
\displaystyle \lefteqn{- \int_\Omega \int_0^T \int_0^1 q_w^+\left( \mu \right) \beta\left( t \right) \d\alpha \d t \psi^\prime\left( x \right) \d x} \next
&\le& \displaystyle \int_0^T \int_\Omega \int_0^1 \eta_w^+\left( \mu \right) \beta^\prime\left( t \right) \psi\left( x \right) \next
&& \displaystyle \vphantom{\int_0^1}
    - \sgn^+\left( \mu - w \right) z^\prime\left( x \right) b\left( \mu \right) \beta\left( t \right) \psi\left( x \right) \d\alpha \d x \d t\next
&\le& \displaystyle C \int_\Omega \psi\left( x \right) \d x
\le -\displaystyle C \int_\Omega x \psi'\left( x \right) \d x
\end{array}
\end{equation}
where $C = C( \mu, z^\prime, b, \beta, w, T )$. But this implies the existence of a set $E^\prime$ of measure zero such that
\begin{equation} \intlabel{monotone function}
x \mapsto \int_0^T \int_0^1 q_w^+\left( \mu \right) \beta\left( t \right) \d\alpha \d t - C x
\quad \left( x \in \Omega \right)
\end{equation}
is non-increasing in $\Omega \setminus E^\prime$. Furthermore there exist a set $E^\pprime$ of measure zero which only depends on $\mu$ and $w$, so that the function in (\intref{monotone function}) is bounded on $\Omega \setminus E^\pprime$. Therefore the function in (\intref{monotone function}) is
monotone and bounded on \ $\Omega \setminus \left(E^\prime \cup E^\pprime \right)$ and the essential limit exists, i.e.
\begin{equation} \intlabel{esslim}
\esslim_{x \uparrow x_r} \int_0^T \int_0^1 q_w^+\left( \mu \right) \beta\left( t \right) \d\alpha \d t
    = \lim_{x \uparrow x_r \above0pt x \notin E^\prime \cup E^\pprime} \int_0^T \int_0^1 q_w^+\left( \mu \right) \beta\left( t \right) \d\alpha \d t.
\end{equation}

Similar as in (\intref{pre-monotony}) we get for $\psi \in C_0^\infty( ] x_l, x_r ] )$:
\begin{equation} \intlabel{pre-esslim}
\begin{array}{rcl}
\lefteqn{\displaystyle -\int_\Omega \int_0^T \int_0^1 q_w^+\left( \mu \right) \beta\left( t \right) \d\alpha \d t \psi^\prime\left( x \right) \d x}\next
&\le& \displaystyle C \int_\Omega \psi\left( x \right) \d x
    + \Lip_{\left[ -C, C \right]}\left( f \right) \int_0^T \eta_w^+\left( u_r\left( t \right) \right) \beta\left( t \right) \d t \psi\left( x_r \right).
\end{array}
\end{equation}

By convolution with a nonnegative kernel we can also apply (\intref{pre-esslim}) to
$\psi_\epsilon( x ) := \frac{1}{\epsilon} ( x - ( x_r - \epsilon ) )^+$, $\epsilon > 0$.
Therefore we obtain
\begin{equation} \intlabel{esslim estimate}
\begin{array}{rcl}
\lefteqn{\displaystyle - \esslim_{x \uparrow x_r} \int_0^T \int_0^1 q_w^+\left( \mu \right) \beta\left( t \right) \d\alpha \d t}\next
&=& \displaystyle
    - \lim_{\epsilon \downarrow 0} \frac{1}{\epsilon} \int_{x_r - \epsilon}^{x_r} \int_0^T \int_0^1 q_w^+\left( \mu \right) \beta\left( t \right) \d\alpha \d t \d x\next
&\le& \displaystyle
    \Lip_{\left[ -C, C \right]}\left( f \right) \int_0^T \eta_w^+\left( u_r\left( t \right) \right) \beta\left( t \right) \d t.
\end{array}
\end{equation}
By further approximation we can also get (\intref{esslim estimate}) for $\beta \in L^1(]0,T[)$, $\beta \ge 0$ a.e. Now let $v \in L^\infty( ] 0, T [ )$, such that $v$ has only a finite number of values $w_i \in \setQ$, i.e.
\begin{equation} \label{elementary-functions}
v = \sum_{i=1}^N w_i \ind_{A_i},
\quad A_i\cap A_j = \emptyset \enspace (i \neq j),
\quad  \bigcup_{i=1,...,N} A_i = ]0,1[.
\end{equation}
Let $\beta \in L^1( ] 0, T [ )$ and $\beta_i := \ind_{A_i} \beta$. Now using  $w_i$ and $\beta_i$ in (\intref{esslim estimate}) and summing over $i$ we obtain:
\begin{equation} \label{pre-result}
\begin{array}{l}
\displaystyle \esslim_{x \uparrow x_r} \int_0^T \int_0^1 \sgn^\pm\left( \mu - v\left( t \right) \right)
   \left( f\left( \mu \right) - f\left( v\left( t \right) \right) \right) \beta\left( t \right) \d\alpha \d t\next
\ge \displaystyle
    - \Lip_{\left[ -C, C \right]}\left( f \right) \int_0^T \left( u_r\left( t \right) - v\left( t \right) \right)^+ \beta\left( t \right) \d t.
\end{array}
\end{equation}
Since any $L^1$-function $v$ is the limit of such elementary functions of type (\ref{elementary-functions}) we obtain (\ref{pre-result}) also for $v \in L^\infty( ] 0, T [ )$.

For the semi-Kru\v zkov entropy pairs $( \eta_w^-, q_w^- )$ and the limits $x \downarrow x_l$ we proceed analogously.
\end{proof}

Now we have to show that the entropy process solution has the correct initial data. Again the proof is similar to the one by Otto without source terms (see \cite{malek}, Chapter 2, Lemma~7.41).

\begin{lemma} \label{initial data}
Let $\mu \in L^\infty( \Omega_T \times ] 0, 1 [ )$, $\mu = \mu( x, t, \alpha )$, be an entropy process solution of (\ref{ibvp}). Then we have:
\begin{displaymath}
\esslim_{t \downarrow 0} \int_\Omega \int_0^1 \left\vert \mu\left( x, t, \alpha \right) - u_0\left( x \right) \right\vert \d\alpha \d x = 0.
\end{displaymath}
\end{lemma}

\begin{proof}
Let $\psi \in C_0^\infty( \Omega )$, choose $( \eta_w^+, q_w^+ )$ and $( \eta_w^-, q_w^- )$ in (\ref{process-sol-kruzkov}) and for $\epsilon > 0$ use
\begin{displaymath}
\varphi_\epsilon\left( x, t \right) := \left[ 1 - \frac{1}{\epsilon} t \right] \ind_{\left[ 0, \epsilon \right]}\left( t \right) \psi\left( x \right)
\end{displaymath}
as test function. Adding those two inequalities we obtain
\begin{displaymath}
\begin{array}{c}
\displaystyle \frac{1}{\epsilon} \int_0^\epsilon \int_\Omega \int_0^1
  - \left\vert \mu - w \right\vert \psi\left( x \right)
  + F\left( \mu, w \right) \left( \epsilon - t \right) \psi^\prime\left( x \right) \d\alpha \d x \d t\next
\displaystyle - \frac{1}{\epsilon} \int_0^\epsilon \int_\Omega \int_0^1
  \sgn\left( \mu - w \right) z^\prime\left( x \right)
    b\left( \mu \right) \left( \epsilon - t \right) \psi\left( x \right) \d\alpha \d x \d t\next
\displaystyle + \int_\Omega \left\vert u_0\left( x \right) - w \right\vert \psi\left( x \right) \d x
\ge 0.
\end{array}
\end{displaymath}
where
\begin{equation} \label{definition-F}
F\left( a, b \right): = \sgn\left( a - b \right) \left( f\left( a \right) - f\left( b \right) \right) \quad \forall a, b \in \setR.
\end{equation}
 Now we take the limes inferior with respect to $\epsilon$ and obtain
\begin{equation} \label{result for w}
\esslimsup_{t \downarrow 0} \int_\Omega \int_0^1 \left\vert \mu - w \right\vert \psi\left( x \right) \d\alpha \d x
\le \int_\Omega \left\vert u_0\left( x \right) - w \right\vert \psi\left( x \right) \d x.
\end{equation}

Similar as in (\ref{pre-result}) we get that (\ref{result for w}) for all $w \in \setR$ implies
\begin{equation} \label{result for v}
\esslimsup_{t \downarrow 0} \int_\Omega \int_0^1 \left\vert \mu - v \right\vert \psi\left( x \right) \d\alpha \d x
\le \int_\Omega \left\vert u_0\left( x \right) - v \right\vert \psi\left( x \right) \d x
\end{equation}
for all $v \in L^\infty( \Omega )$ and all $\psi \in L^1( \Omega )$, $\psi \ge 0$ a.e.

Choosing $u_0$ for $v$ in (\ref{result for v}) and using $\psi = 1$ as test function we get the statement of the lemma.
\end{proof}

The proofs of the following two lemmata are similar to those in \cite{vovelle} and therefore we omit them. They are mainly based on \lemmaname~\ref{trace}.

\begin{lemma}
Let $\mu \in L^\infty( \Omega_T \times ] 0, 1 [ )$ be an entropy  process solution of (\ref{ibvp}). Using $F$ given by (\ref{definition-F}) we have for all $\varphi \in C_0^\infty( \overline\Omega \times ] 0, T [ )$, $\varphi \ge 0$:
\begin{displaymath}
\begin{array}{c}
\displaystyle \int_0^T \int_\Omega \int_0^1 \left\vert \mu - k \right\vert \partial_t \varphi + F\left( \mu, k \right) \partial_x \varphi \d\alpha \d x \d t\next
\displaystyle - \int_0^T \int_\Omega \int_0^1 \sgn\left( \mu - k \right) z^\prime\left( x \right) b\left( \mu \right) \varphi \d\alpha \d x \d t\next
\displaystyle
  + \esslim_{x \uparrow x_r} \int_0^T \int_0^1 F\left( \mu, u_r\left( t \right) \right) \varphi\left( x_r, t \right) \d\alpha \d t\next
\displaystyle
  - \esslim_{x \downarrow x_l} \int_0^T \int_0^1 F\left( \mu, u_l\left( t \right) \right) \varphi\left( x_l, t \right) \d\alpha \d t\next
\displaystyle
  - \int_0^T F\left( u_r\left( t \right), k \right) \, \varphi\left( x_r, t \right) \, dt
    + \int_0^T F\left( u_l\left( t \right), k \right) \, \varphi\left( x_l, t \right) \, dt
\ge 0,
\end{array}
\end{displaymath}
where $\mu = \mu( x, t, \alpha )$ and $\varphi = \varphi( x, t )$.
\end{lemma}

\begin{lemma}
Let $\mu, \nu \in L^\infty( \Omega_T \times ] 0, 1 [ )$ be two entropy process solutions of (\ref{ibvp}) with respect to the initial data $u_0 \in L^\infty( \Omega )$ and the boundary data $u_l, u_r \in L^\infty( ] 0, T [ )$. Then donoting $\mu = \mu( x, t, \alpha )$, $\nu = \nu( x, t, \beta )$ we have for all nonnegative $\varphi \in C_0^\infty( \overline\Omega \times ] 0, T [ )$, $\varphi = \varphi( x, t )$:
\begin{equation} \label{uniqueness-inequality}
\begin{array}{c}
\displaystyle \int_{\Omega_T} \int_0^1 \int_0^1 \left\vert \mu - \nu \right\vert \partial_t \varphi
    + F\left( \mu, \nu \right) \partial_x \varphi \d\beta \d\alpha \d(x,t)\next
\displaystyle - \int_{\Omega_T} \int_0^1 \int_0^1 \sgn\left( \mu - \nu \right) z^\prime\left( x \right)
      \left( b\left( \mu \right) - b\left( \nu \right) \right) \varphi \d\beta \d\alpha \d(x,t)
\ge 0,
\end{array}
\end{equation}
where $F$ is given by  (\ref{definition-F}).
\end{lemma}

Now we are ready to prove the uniqueness of the entropy process solution. The proof without the source term can be found in \cite{vovelle}, Theorem~2. For controlling the source term we use similar ideas as in \cite{bardos}, Theorem~5.

\begin{numproof}[Proof of \theoremname~\ref{uniqueness}]
Let $\psi \in C_0^\infty( ] 0, T [ )$, $\psi \ge 0$. In (\ref{uniqueness-inequality}) we choose
$( x, t ) \mapsto \psi( t )$ as a test function and obtain with $\mu = \mu( x, t, \alpha )$, $\nu = \nu( x, t, \beta )$:
\begin{equation} \intlabel{term 1}
\begin{array}{c}
\displaystyle \int_{\Omega_T} \int_0^1 \int_0^1
    \left\vert \mu - \nu \right\vert \partial_t \psi\left( t \right) \d\beta \d\alpha \d(x,t)\next
\displaystyle - \int_{\Omega_T} \int_0^1 \int_0^1
    \sgn\left( \mu - \nu \right) z^\prime\left( x \right)
    \left( b\left( \mu \right) - b\left( \nu \right) \right) \psi\left( t \right) \d\beta \d\alpha \d (x,t)
\ge 0.
\end{array}
\end{equation}

Using
\begin{displaymath}
g\left( t \right) :=
    \int_\Omega \int_0^1 \int_0^1 \left\vert \mu\left( x, t, \alpha \right) - \nu\left( x, t, \beta \right) \right\vert \d\beta \d\alpha \d x
\quad \mbox{for} \enspace t \in \left] 0, T \right[,
\end{displaymath}
(\intref{term 1}) can be written as
\begin{displaymath}
- \int_0^T g\left( t \right) \partial_t \psi\left( t \right) \d t
\le \left\Vert z^\prime \right\Vert_{L^\infty\left( \Omega \right)} \left\Vert b^\prime \right\Vert_{L^\infty\left( \setR \right)}
    \int_0^T g\left( t \right) \psi\left( t \right) \d t.
\end{displaymath}
Integration by parts on the right hand side implies that there exists a set $E \subset ] 0, T [$ of measure zero, such that
\begin{displaymath}
t \mapsto g\left( t \right)
    - \left\Vert z^\prime \right\Vert_{L^\infty\left( \Omega \right)} \left\Vert b^\prime \right\Vert_{L^\infty\left( \setR \right)} \int_0^t g\left( \xi \right) \d\xi
\end{displaymath}
is non-increasing on $] 0, T [ \setminus E$, i.e. for all $t_1, t_2 \in ] 0, T [ \setminus E$, $t_1 < t_2$, we have:
\begin{equation} \intlabel{term 2}
g\left( t_2 \right) \le g\left( t_1 \right) + \left\Vert z^\prime \right\Vert_{L^\infty\left( \Omega \right)}
    \left\Vert b^\prime \right\Vert_{L^\infty\left( \setR \right)} \int_{t_1}^{t_2} g\left( \xi \right) \d\xi.
\end{equation}
Now we can apply the Gronwall Lemma  and (\intref{term 2}) implies
\begin{equation} \intlabel{term 3}
\begin{array}{rcl}
\displaystyle \lefteqn{\left\Vert \mu\left( \cdot, t_2, \cdot \right)
    - \nu\left( \cdot, t_2, \cdot \right) \right\Vert_{L^1\left( \Omega \times \left] 0, 1 \right[^2 \right)}}\next
&\le& \displaystyle \left\Vert \mu\left( \cdot, t_1, \cdot \right) - \nu\left( \cdot, t_1, \cdot \right) \right\Vert_{L^1\left( \Omega \times \left] 0, 1 \right[^2 \right)}
    e^{\left\Vert z^\prime \right\Vert_{L^\infty\left( \Omega \right)} \left\Vert b^\prime \right\Vert_{L^\infty\left( \setR \right)} \left( t_2 - t_1 \right)}.
\end{array}
\end{equation}
Since the entropy process solution respects the initial data in the sense of \lemmaname~\ref{initial data} we have:
\begin{displaymath}
\begin{array}{rcl}
\displaystyle \lefteqn{\esslimsup_{t_1 \downarrow 0} \left\Vert \mu\left( \cdot, t_1, \cdot \right) - \nu\left( \cdot, t_1, \cdot \right)
    \right\Vert_{L^1\left( \Omega \times \left] 0, 1 \right[^2 \right)}}\next
&\le& \displaystyle \esslim_{t_1 \downarrow 0}
  \left[ \left\Vert \mu\left( \cdot, t_1, \cdot \right) - u_0 \right\Vert_{L^1\left( \Omega \times \left] 0, 1 \right[ \right)}
  + \left\Vert \nu\left( \cdot, t_1, \cdot \right) - u_0 \right\Vert_{L^1\left( \Omega \times \left] 0, 1 \right[ \right)} \right]\next
&=& \displaystyle 0,
\end{array}
\end{displaymath}
and therefore from (\intref{term 3}) for almost all $t_2 \in \left] 0, T \right[$:
\begin{displaymath}
\left\Vert \mu\left( \cdot, t_2, \cdot \right) - \nu\left( \cdot, t_2, \cdot \right) \right\Vert_{L^1\left( \Omega \times \left] 0, 1 \right[^2 \right)} \le 0.
\end{displaymath}

But this proves  $\mu\left( x, t, \alpha \right) = \nu\left( x, t, \beta \right)$ for almost all  $\left( x, t \right) \in \Omega_T$ and almost all $\alpha, \beta \in \left] 0, 1 \right[.$ Therefore $\mu$ and $\nu$ do not depend on $\alpha$ and $\beta$ respectively. Setting
\begin{displaymath}
u\left( x, t \right) := \int_0^1 \mu\left( x , t, \alpha \right) \d\alpha
\end{displaymath}
we get for almost all $( x, t, \alpha ) \in \Omega_T \times ] 0, 1 [$:
\begin{displaymath}
\mu\left( x, t, \alpha \right) =  u\left( x, t \right) = \nu\left( x, t, \alpha \right).
\end{displaymath}
Then, by \definitionname s~\ref{def-entropy-sol} and \ref{def-process-sol}, we see that $u$ is an entropy solution of (\ref{ibvp}).
\end{numproof}

\section{The Well-Balanced Scheme} \label{well}

In this section we will describe the well-balanced scheme for the initial boundary value problem (\ref{ibvp}). It  was originally developed by Greenberg et al.\ in \cite{greenberg} for the initial value problem. For the Engquist-Osher numerical flux Perthame et al. (see \cite{perthame}) could prove convergence of the numerical solution to the entropy solution. Here we are going to generalize this result to the initial boundary value problem (\ref{ibvp}).

First let us fix the notation for the discretization. Let $\Delta x > 0$ such that $\frac{\vert \Omega \vert}{\Delta x} \in \setN$, $\Omega = ] x_l, x_r [$, $j_l \in \setZ$ and $x_j := x_l + (j - j_l+\frac{1}{2}) \Delta x$, $C_j := ] x_{j-\frac{1}{2}}, x_{j+\frac{1}{2}} [$, choose  $J \subset \setZ$ such
that $\overline\Omega = \bigcup_{j \in J} \overline C_j$ and define  $j_r := \max J$.

\begin{center}
\setlength{\unitlength}{0.7\unitlength}
\begin{picture}(420,50)
\put(10,20){\line(1,0){90}} \put(160,20){\line(1,0){100}}
\put(320,20){\line(1,0){90}}
\put(100,10){\makebox(60,20)[c]{$\cdots$}}
\put(260,10){\makebox(60,20)[c]{$\cdots$}}

\put(10,16){\line(0,1){8}} \put(90,16){\line(0,1){8}}
\put(170,16){\line(0,1){8}} \put(250,16){\line(0,1){8}}
\put(330,16){\line(0,1){8}} \put(410,16){\line(0,1){8}}
\put(50,18){\line(0,1){4}} \put(210,18){\line(0,1){4}}
\put(370,18){\line(0,1){4}}

\put(0,0){\makebox(20,14)[t]{$x_l$}}
\put(400,0){\makebox(20,14)[t]{$x_r$}}

\put(40,0){\makebox(20,14)[t]{$x_{j_l}$}}
\put(200,0){\makebox(20,14)[t]{$x_j$}}
\put(360,0){\makebox(20,14)[t]{$x_{j_r}$}}

\put(160,25){\makebox(100,25)[b]{$\overbrace{\makebox[80\unitlength]{}}^{C_j}$}}
\end{picture}
\end{center}

Let $\Delta t > 0$ such that $N_T := \frac{T}{\Delta t} \in \setN$ and $t^n := n \Delta t$. For any function $\varphi \in C^0( \setR^2 )$ we define $\varphi^n_j := \varphi\left( x_j, t^n \right)$ and the piecewise constant function
\begin{displaymath}
\overline\varphi\left( x, t \right) := \varphi^n_j \quad t \in [ t^n, t^{n+1} ), \enspace x \in C_j.
\end{displaymath}

The main idea for the well-balanced schemes consists in the following fact. If $v$ is a stationary solution of (\ref{ibvp}) then $v$ satifies
\begin{displaymath}
\partial_x \left( D\left( v\left( x \right) \right) + z\left( x \right) \right) = 0 \quad \mbox{for a.a.} \enspace x \in \setR,
\end{displaymath}
(for the definition of $D$ see (\ref{definition D})) which is equivalent to
\begin{equation} \label{equilibrium-D}
D\left( v \right) + z\left( x \right) = c \in \setR \quad \mbox{for all} \enspace x \in \setR.
\end{equation}
This property is the main building block for the numerical scheme.

The standard form of a numerical scheme in conservation form for the partial differential equation in (\ref{ibvp}) is
\begin{equation} \label{standard-scheme}
u_j^{n+1} := u_j^n - \frac{\Delta t}{\Delta x} \left( g\left( u_j^n, u_{j+1}^n \right)- g\left( u_{j-1}^n, u_j^n \right) \right) - \Delta t z^\prime_j u_j^n
\quad \mbox{for} \enspace j\in J
\end{equation}
where $g$ is a numerical flux and $z^\prime_j$ a discretisation of $z^\prime$ in $C_j$, e.g.\ the average of $z^\prime$ on the cell.

This scheme is very inefficient, especially for the approximation of stationary solutions. The well-balanced schemes are much better and are defined as follows.

\begin{definition}[Well-Balanced Scheme] \label{def-scheme}
Let
\begin{equation} \label{ini}
u^0_j := \frac{1}{\Delta x} \int_{C_j} u_0\left( x \right) \d x
\end{equation}
and define $u^n_j$ for $j \notin J$ by
\begin{equation} \label{ex-ini}
u^n_j \:= \left\lbrace
\begin{array}{ll}
\displaystyle u^n_l = \frac{1}{\Delta t} \int_{t^n}^{t^{n+1}} u_l\left( t \right) \d t & \displaystyle \enspace \mbox{for} \enspace j < j_l,\next
\displaystyle u^n_r = \frac{1}{\Delta t} \int_{t^n}^{t^{n+1}} u_r\left( t \right) \d t & \displaystyle \enspace \mbox{for} \enspace j > j_r.
\end{array}
\right.
\end{equation}
Similarly let $z_j$ be given by
\begin{displaymath}
z_j := \left\lbrace
\begin{array}{ll}
\displaystyle z_{j_l} & \displaystyle \enspace \mbox{for} \enspace j < j_l,\next
\displaystyle {1 \over \Delta x} \int_{C_j} z\left( x \right) \, dx & \displaystyle \enspace \mbox{for} \enspace j \in J,\next
\displaystyle z_{j_r} & \displaystyle \enspace \mbox{for} \enspace j > j_r.
\end{array}\right.
\end{displaymath}

Now assume that $( u^n_j )_{j \in \setZ}$ is already defined. Then the values $u^{n+1}_j$ for the new time step of the well-balanced scheme are given by
\begin{equation} \label{scheme}
u^{n+1}_j: = u^n_j - \frac{\Delta t}{\Delta x} \left( g\left( u^n_j, u^n_{j+1,-} \right) - g\left( u^n_{j-1,+}, u^n_j \right) \right) \quad \mbox{for} \enspace j\in J
\end{equation}
where, due to (\ref{equilibrium-D}), $u^n_{j+1,-}$ and $u^n_{j-1,+}$ are defined by
\begin{equation} \label{discrete equilibrium-states}
\begin{array}{rcl}
D\left( u_{j-1,+} \right) + z_j &=& D\left( u_{j-1} \right) + z_{j-1},\next
D\left( u_{j+1,-} \right) + z_j &=& D\left( u_{j+1} \right) + z_{j+1}.
\end{array}
\end{equation}

Using the discrete data  $( u^n_j )_{j \in \setZ, 0 \le n < N_T}$ we define the numerical solution $u_{\Delta x} \in L^\infty( \setR \times [ 0, T [ )$ by
\begin{equation}
u_{\Delta x}\left( x, t \right) := u^n_j \quad \mbox{for} \enspace t \in [ t^n, t^{n+1} [, \enspace n < N_T, \enspace x \in C_j, \enspace j \in \setZ.
\end{equation}
\end{definition}

\begin{remark}
Under condition (\ref{assumptions-D}) there always exist unique solutions $u_{j-1,+}$ and $ \,u_{j+1,-}$ of (\ref{discrete equilibrium-states}).
\end{remark}

\begin{remark}
If the values $u^n_j$ are ``stationary'', i.e.\ $D( u^n_j ) + z_j = D( u^n_{j-1} ) + z_{j-1}$ for all $j \in J$, then we obtain
\begin{displaymath}
u^n_{j-1,+} = u^n_j, \quad u^n_{j+1,-} = u^n_j
\end{displaymath}
as solutions of (\ref{discrete equilibrium-states}) and therefore $u^{n+1}_j = u^n_j$ for all $j \in J$.
\end{remark}

Now we assume that the numerical flux is given by the Engquist-Osher flux
\begin{equation} \label{engquist-osher flux}
g\left( u, v \right) = \int_0^u {f^\prime}^+\left( \xi \right) \d\xi - \int_0^v {f^\prime}^-\left( \xi \right) \d\xi + f\left( 0 \right)
\quad \mbox{for} \enspace u, v \in \setR.
\end{equation}
and set
\begin{equation} \label{engquist-osher scheme}
H\left( u, v, w \right) := v - \frac{\Delta t}{\Delta x} \left( g\left( v, w \right) - g\left( u, v \right) \right)
\quad \mbox{for}\enspace u, v, w \in \setR.
\end{equation}

Before proving the boundedness of the numerical solutions we will prove that the numerical values obtained by the Engquist-Osher scheme remain in the convex hull of its data.

As in \cite{perthame} we use the following ``density-function'' in order to move the integration bounds to the integrand:
\begin{displaymath}
\chi_s\left( \xi \right) := \left\lbrace
\begin{array}{ll}
1 & \enspace \mbox{for} \enspace 0 < \xi < s,\\
-1 & \enspace \mbox{for} \enspace s < \xi < 0,\\
0 & \enspace \mbox{for} \enspace otherwise.
\end{array}\right.
\end{displaymath}

Now we can write the Engquist-Osher scheme (\ref{engquist-osher scheme}) as an integral by defining
\begin{equation} \label{definition-h}
\begin{array}{rcl}
\displaystyle h\left( u, v, w \right)\left( \xi \right)
&:=& \displaystyle \chi_v\left( \xi \right) - \frac{\Delta t}{\Delta x}
    \left( {f^\prime}^+\left( \xi \right) \chi_v\left( \xi \right)
         - {f^\prime}^-\left( \xi \right) \chi_w\left( \xi \right) \right)\next
&& \displaystyle + \frac{\Delta t}{\Delta x}
    \left( {f^\prime}^+\left( \xi \right) \chi_u\left( \xi \right)
         + {f^\prime}^-\left( \xi \right) \chi_v\left( \xi \right) \right)
\end{array}
\end{equation}
for all $u, v, w, \xi \in \setR$. Then (\ref{engquist-osher scheme}) can be written as
\begin{displaymath}
H\left( u, v, w \right) = \int_\setR h\left( u, v, w \right)\left( \xi \right) \d\xi
\quad \mbox{for} \enspace u, v, w \in \setR.
\end{displaymath}

This form can be used to show that the values produced by the Engquist-Osher scheme are the convex hull of its arguments.
\begin{lemma} \label{properties-h}
Let  $h : \setR^3 \to L^1\left( \setR \right)$ be defined by (\ref{definition-h}), $K = \left[ -C, C \right]$, $\,u, v, w \in K$ and assume that the
CFL-condition
\begin{displaymath}
\left\Vert f^\prime \right\Vert_{L^\infty\left( K \right)} \frac{\Delta t}{\Delta x} \le 1
\end{displaymath}
holds. Then for all  $\xi \in \setR$ we have
\begin{eqnarray*}
&h\left( u, v, w \right)\left( \xi \right) \in \conv\left\lbrace \chi_u\left( \xi \right), \chi_v\left( \xi \right), \chi_w\left( \xi \right) \right\rbrace,\\
&0 \le \sgn\left( \xi \right) h\left( u, v, w \right)\left( \xi \right) \le 1
\end{eqnarray*}
and consequently
\begin{equation} \label{max-H}
\left\vert H\left( u, v, w \right) \right\vert \le \max\left\lbrace \left\vert u \right\vert, \left\vert v \right\vert, \left\vert w \right\vert \right\rbrace.
\end{equation}
\end{lemma}

\begin{proof}
The proof uses only standard means, see also \cite{perthame}.
\end{proof}

In the next step we will prove the boundedness of the numerical solutions.

\begin{lemma}[Boundedness of the Numerical Solutions] \label{lemma-boundedness}
Let (\ref{ass}), (\ref{assumptions-D}), (\ref{engquist-osher flux}) be satisfied,
\begin{eqnarray}
M &:=& \max\left\lbrace \left\Vert u_0 \right\Vert_{L^\infty\left( \Omega \right)},
    \left\Vert u_l \right\Vert_{L^\infty\left( \left] 0, T \right[ \right)},
    \left\Vert u_r \right\Vert_{L^\infty\left( \left] 0, T \right[ \right)} \right\rbrace,\\
\label{L-infty-bound}
C^{\Delta x}_T &:=& M e^{2 T \left\Vert b^\prime \right\Vert_{L^\infty\left( \setR \right)} \left\Vert z^\prime \right\Vert_{L^\infty\left( \Omega \right)}}\nonumber\\
&&+  \Delta x \frac{\left\Vert z^\prime \right\Vert_{L^\infty\left( \Omega \right)}}{\inf_\setR D^\prime}
    e^{4 T \left\Vert b^\prime \right\Vert_{L^\infty\left( \setR \right)} \left\Vert z^\prime \right\Vert_{L^\infty\left( \Omega \right)}}\\
&&+ \left\vert b\left( 0 \right) \right\vert
    e^{2 T \left( \left\Vert b^\prime \right\Vert_{L^\infty\left( \setR \right)} + 1 \right)
    \left\Vert z^\prime \right\Vert_{L^\infty\left( \Omega \right)}},\nonumber\\
K^{\Delta x}_T &:=& \left[ -C^{\Delta x}_T, C^{\Delta x}_T \right]\nonumber
\end{eqnarray}
and we assume the CFL-condition $\Lip_{K^{\Delta x}_T}\left( f \right)\frac{\Delta t}{\Delta x} \le 1.$  Then we obtain for all $n \, \Delta t \le T$:
\begin{equation} \intlabel{estimate-u}
\sup_{j \in J}\vert u^n_j \vert \le C^{\Delta x}_T.
\end{equation}
Furthermore we have for all $\left( n+1 \right) \Delta t \le T$:
\begin{equation} \intlabel{estimate-equilibrium-states}
\sup_{j \in J}\vert u^n_{j-1,+} \vert \le C^{\Delta x}_T,
\quad \sup_{j \in J}\vert u^n_{j+1,-} \vert \le C^{\Delta x}_T.
\end{equation}
\end{lemma}

\begin{proof}
Let  $u, v, w \in \setR$ be arbitrary. Then we have with $b_0 := b( 0 )$:
\begin{eqnarray*}
\left\vert g\left( u, v \right) - g\left( u, w \right) \right\vert
&=& \left\vert \int_v^w {f^\prime}^-\left( s \right) \d s \right\vert
\le \int_{\left[ v, w \right]} \left\vert {f^\prime}^-\left( s \right) \right\vert \d s\\
&\le& \int_{\left[ v, w \right]} \left\vert f^\prime\left( s \right) \right\vert \d s
= \int_{\left[ v, w \right]} \left\vert D^\prime\left( s \right) \, b\left( s \right) \right\vert \d s\\
&\le& \left\Vert b^\prime \right\Vert_\infty \int_{\left[ v, w \right]} D^\prime\left( s \right) \left\vert s \right\vert \d s
    + \left\vert b_0 \right\vert \int_{\left[ v, w \right]} D^\prime\left( s \right) \d s\\
&\le& \left[ \left\Vert b^\prime \right\Vert_\infty \max\left\lbrace \left\vert v \right\vert, \left\vert w \right\vert \right\rbrace
    + \left\vert b_0 \right\vert \right] \int_{\left[ v, w \right]} D^\prime\left( s \right)\d s\\
&\le& \left[ \left\Vert b^\prime \right\Vert_\infty \max\left\lbrace \left\vert v \right\vert, \left\vert w \right\vert \right\rbrace
     + \left\vert b_0 \right\vert \right] \left\vert D\left( v \right) - D\left( w \right) \right\vert.
\end{eqnarray*}
The mean value theorem implies:
\begin{equation} \intlabel{diff}
\begin{array}{rcl}
\displaystyle \left\vert u^n_{j+1} - u^n_{j+1,-} \right\vert
&\le& \displaystyle \frac{\vert D(u^n_{j+1}) - D(u^n_{j+1,-}) \vert}{\inf_\setR D^\prime}\next
&=& \displaystyle \frac{\left\vert z_j - z_{j+1} \right\vert}{\inf_\setR D^\prime}
\le \displaystyle \Delta x \frac{\left\Vert z^\prime \right\Vert_\infty}{\inf_\setR D^\prime}
\end{array}
\end{equation}
and therefore:
\begin{eqnarray*}
\lefteqn{\left\vert g\left( u^n_j, u^n_{j+1,-} \right) - g\left( u^n_j, u^n_{j+1} \right) \right\vert}\\
&\le& \left[ \left\Vert b^\prime \right\Vert_\infty \left( \left\vert u^n_{j+1} \right\vert
    + \Delta x \frac{\left\Vert z^\prime \right\Vert_\infty}{\inf_\setR D^\prime} \right)  + \left\vert b_0 \right\vert \right]
    \Delta x \left\Vert z^\prime \right\Vert_\infty.
\end{eqnarray*}

In a similar way we obtain:
\begin{eqnarray*}
\lefteqn{\left\vert g\left( u^n_{j-1,+}, u^n_j \right) - g\left( u^n_{j-1}, u^n_j \right) \right\vert}\\
&\le& \left[ \left\Vert b^\prime \right\Vert_\infty \left( \left\vert u^n_{j-1} \right\vert
    + \Delta x \frac{\left\Vert z^\prime \right\Vert_\infty}{\inf_\setR D^\prime} \right) + \left\vert b_0 \right\vert \right]
    \Delta x \left\Vert z^\prime \right\Vert_\infty.
\end{eqnarray*}

Now by induction with respect to $n \in \setN$ such that  $n \Delta t \le T$ we will show:
\begin{equation} \intlabel{goal}
\begin{array}{rcl}
\displaystyle
\sup_{j \in J} \left\vert u^n_j \right\vert
&\le& \displaystyle
    M e^{2 n \Delta t \left\Vert b^\prime \right\Vert_\infty \left\Vert z^\prime \right\Vert_\infty}
    +  \Delta x \frac{\left\Vert z^\prime \right\Vert_\infty}{\inf_\setR D^\prime}
        e^{4 n \Delta t \left\Vert b^\prime \right\Vert_\infty \left\Vert z^\prime \right\Vert_\infty}\next
&&\displaystyle
    + \left\vert b_0 \right\vert
        e^{2 n \Delta t \left( \left\Vert b^\prime \right\Vert_\infty + 1 \right)
        \left\Vert z^\prime \right\Vert_\infty},
\end{array}
\end{equation}
which will prove the statement (\intref{estimate-u}) of the lemma. Obviously (\intref{goal}) holds for the initial data, i.e. $n = 0$.

Assume that (\intref{goal}) holds for  $n \in \setN$, $\left( n+1 \right) \Delta t \le T$. Since the maximum principle holds for the Engquist-Osher flux $g(u,v)$ (see
(\ref{max-H})), the CFL-condition holds and $\Lip_{K^{\Delta x}_T}\left( g \right) \le \Lip_{K^{\Delta x}_T}\left( f \right)$ we obtain for all  $j \in J$:
\begin{eqnarray*}
\left\vert u^{n+1}_j \right\vert
&\le& \left\vert u^n_j - {\Delta t \over \Delta x} \left( g\left( u^n_j, u^n_{j+1} \right) - g\left( u^n_{j-1}, u^n_j \right) \right) \right\vert\\
&&+ \left\vert {\Delta t \over \Delta x} \left( g\left( u^n_j, u^n_{j+1,-} \right) - g\left( u^n_j, u^n_{j+1} \right) \right) \right\vert\\
&&+ \left\vert {\Delta t \over \Delta x} \left( g\left( u^n_{j-1,+}, u^n_j \right) - g\left( u^n_{j-1}, u^n_j  \right) \right) \right\vert\\
&\le& \sup_{j \in \setZ} \left\vert u^n_j \right\vert
    + 2 \Delta t \left[ \left\Vert b^\prime \right\Vert_\infty \left[ \sup_{j \in \setZ} \left\vert u^n_j \right\vert
        + \frac{\Delta x \left\Vert z^\prime \right\Vert_\infty}{\inf_\setR D^\prime} \right]
        + \left\vert b_0 \right\vert \right] \left\Vert z^\prime \right\Vert_\infty\\
&\le& \max\left\lbrace \sup_{j \in J} \left\vert u^n_j \right\vert, \left\vert u^n_l \right\vert, \left\vert u^n_r \right\vert \right\rbrace
    \left( 1 + 2 \Delta t \left\Vert b^\prime \right\Vert_\infty \left\Vert z^\prime \right\Vert_\infty \right)\\
&&+ 2 \Delta t \left\vert b_0 \right\vert \left\Vert z^\prime \right\Vert_\infty
    + 2 \Delta x \Delta t \left\Vert b^\prime \right\Vert_\infty
        \frac{\left\Vert z^\prime \right\Vert^2_\infty}{\inf_\setR D^\prime}.
\end{eqnarray*}

If
\begin{math}
\max\lbrace \sup_{j \in J} \vert u^n_j \vert, \vert u^n_l \vert, \vert u^n_r \vert \rbrace
= \vert u^n_l \vert,
\end{math}
we have
\begin{eqnarray*}
\left\vert u^{n+1}_j \right\vert
&\le& \left\Vert u_l \right\Vert_{L^\infty\left( \left] 0, T \right[ \right)}
    e^{2 \Delta t \left\Vert b^\prime \right\Vert_\infty \left\Vert z^\prime \right\Vert_\infty}
    \underbrace{e^{2 n \Delta t \left\Vert b^\prime \right\Vert_\infty
        \left\Vert z^\prime \right\Vert_\infty}}_{\ge 1}\\
&&+ \frac{\Delta x \left\Vert z^\prime \right\Vert_\infty}{\inf_\setR D^\prime}
    e^{4 \Delta t \left\Vert b^\prime \right\Vert_\infty \left\Vert z^\prime \right\Vert_\infty}
    \underbrace{e^{4 n \Delta t \left\Vert b^\prime \right\Vert_\infty
        \left\Vert z^\prime \right\Vert_\infty}}_{\ge 1}\\
&&+ \left\vert b_0 \right\vert e^{2 \Delta t \left\Vert z^\prime \right\Vert_\infty}
    e^{2 n \Delta t \left( \left\Vert b^\prime \right\Vert_\infty + 1 \right)
    \left\Vert z^\prime \right\Vert_\infty}\\
&\le& M e^{2 \left( n+1 \right) \Delta t \left\Vert b^\prime \right\Vert_\infty \left\Vert z^\prime \right\Vert_\infty}
 +  \frac{\Delta x \left\Vert z^\prime \right\Vert_\infty}{\inf_\setR D^\prime}
        e^{4 \left( n+1 \right) \Delta t \left\Vert b^\prime \right\Vert_\infty \left\Vert z^\prime \right\Vert_\infty}\\
&&+ \left\vert b_0 \right\vert
    e^{2 \left( n+1 \right) \Delta t \left( \left\Vert b^\prime \right\Vert_\infty + 1 \right)
    \left\Vert z^\prime \right\Vert_\infty}.
\end{eqnarray*}

In a similar way the case
\begin{math}
\max\lbrace \sup_{j \in J} \vert u^n_j \vert, \vert u^n_l \vert, \vert u^n_r \vert \rbrace
= \vert u^n_r \vert
\end{math}
can be handled.

If
\begin{math}
\max\lbrace \sup_{j \in J} \vert u^n_j \vert, \vert u^n_l \vert, \vert u^n_r \vert \rbrace
= \sup_{j \in J} \vert u^n_j \vert
\end{math}
we use (\intref{goal}) to obtain:
\begin{eqnarray*}
\left\vert u^{n+1}_j \right\vert
&\le& M  e^{2 n \Delta t \left\Vert b^\prime \right\Vert_\infty \left\Vert z^\prime \right\Vert_\infty}
    \left( 1 + 2 \Delta t \left\Vert b^\prime \right\Vert_\infty \left\Vert z^\prime \right\Vert_\infty \right)\\
&&+ \frac{\Delta x \left\Vert z^\prime \right\Vert_\infty}{\inf_\setR D^\prime}
    e^{4 n \Delta t \left\Vert b^\prime \right\Vert_\infty \left\Vert z^\prime \right\Vert_\infty}
        \left( 1 + 2 \Delta t \left\Vert b^\prime \right\Vert_\infty \left\Vert z^\prime \right\Vert_\infty \right)\\
&&+ \left\vert b_0 \right\vert
    e^{2 n \Delta t \left( \left\Vert b^\prime \right\Vert_\infty + 1 \right) \left\Vert z^\prime \right\Vert_\infty}
    \left( 1 + 2 \Delta t \left\Vert b^\prime \right\Vert_\infty \left\Vert z^\prime \right\Vert_\infty \right)\\
&&+ 2 \Delta t \left\vert b_0 \right\vert \left\Vert z^\prime \right\Vert_\infty
  + 2 \Delta x \Delta t \left\Vert b^\prime \right\Vert_\infty
    \frac{\left\Vert z^\prime \right\Vert^2_\infty}{\inf_\setR D^\prime}.\\
&\le& M e^{2 \left( n+1 \right) \Delta t \left\Vert b^\prime \right\Vert_\infty \left\Vert z^\prime \right\Vert_\infty}
    + \frac{\Delta x \left\Vert z^\prime \right\Vert_\infty}{\inf_\setR D^\prime}
        e^{4 \left( n+1 \right) \Delta t \left\Vert b^\prime \right\Vert_\infty \left\Vert z^\prime \right\Vert_\infty}\\
&&+ \left\vert b_0 \right\vert
    e^{2 \left( n+1 \right) \Delta t \left( \left\Vert b^\prime \right\Vert_\infty + 1 \right) \left\Vert z^\prime \right\Vert_\infty}.
\end{eqnarray*}
Therefore we have proved (\intref{goal}) for $n \Delta t \le T$. The estimates in (\intref{estimate-equilibrium-states}) now follow from (\intref{diff}), (\intref{goal}) and
$\left( n+1 \right) \Delta t \le T$.
\end{proof}

Notice that the second term in (\ref{L-infty-bound}) is necessary, although it is missing in the proof for the initial value problem in \cite{perthame}, Lemma~3.1. This is shown in the diplom thesis \cite{nolte}.

A main tool for the convergence proof in Section \ref{main} will be a cell entropy inequality. This means for an arbitrary entropy pair $\left( \eta, q \right)$ we have to find a numerical entropy flux $G$, such that for $C > 0$ and for $\frac{\Delta t}{\Delta x}$ sufficiently small
\begin{equation} \label{cell-entropy-inequality}
\eta\left( H\left( u, v, w \right) \right) - \eta\left( v \right) + \frac{\Delta t}{\Delta x} \left( G\left( v, w \right) - G\left( u, v \right) \right)
\le 0
\end{equation}
holds for all $u, v, w \in \left[ -C, C \right]$.

To prove the inequality in (\ref{cell-entropy-inequality}) we will need a lemma by Brenier (see \cite{brenier}).

\begin{lemma}[Brenier \cite{brenier}] \label{lemma-brenier}
Let $f \in L^1\left( \setR \right)$ such that $0 \le \sgn\left( \xi \right) f\left( \xi \right) \le 1$ for all $\xi \in \setR$ and let $h \in C^2\left( \setR \right)$ be convex and Lipschitz-continuous. Then we have:
\begin{displaymath}
h\left( \int_\setR f\left( \xi \right) \d\xi \right) - h\left( 0 \right)
\le \int_\setR h^\prime\left( \xi \right) f\left( \xi \right) \d\xi.
\end{displaymath}
\end{lemma}

Now we can prove the cell entropy inequality (\ref{cell-entropy-inequality}), which generalizes a result due to Perthame et al., see \cite{perthame}.

\begin{lemma} \label{lemma-cell-entropy-inequality}
Let  $\left\vert u \right\vert, \left\vert v \right\vert, \left\vert w \right\vert \le C \in \setR$, $K := \left[ -C, C \right]$, $\left( \eta, q \right)$ be an entropy pair and assume the CFL-condition
\begin{equation} \intlabel{cfl}
\left\Vert f^\prime \right\Vert_{L^\infty\left( K \right)} \frac{\Delta t}{\Delta x} \le 1.
\end{equation}

Then $G : \setR \times \setR \to \setR$ given by
\begin{equation} \label{definition-G}
G\left( u, v \right) := \int_0^u \eta^\prime\left( \xi \right) {f^\prime}^+\left( \xi \right) \d\xi
    - \int_0^v \eta^\prime\left( \xi \right) {f^\prime}^-\left( \xi \right) \d\xi + q\left( 0 \right)
\end{equation}
is a consistent numerical entropy flux, such that the cell entropy inequality (\ref{cell-entropy-inequality}) is satisfied.
\end{lemma}

\begin{proof}
It is easy to check that $G$ is a consistent numerical entropy flux. Since $\left\vert u \right\vert, \left\vert v \right\vert, \left\vert w \right\vert \le C$ and due to (\intref{cfl}) $\left\vert H\left( u, v, w \right) \right\vert \le C$, we can assume without loss of generality that $\eta$ is Lipschitz-continuous.

Due to the CFL-condition (\intref{cfl}) and \lemmaname~\ref{properties-h} we have for all $\xi \in \setR$ that
\begin{displaymath}
0 \le \sgn\left( \xi \right) h\left( u, v, w \right)\left( \xi \right) \le 1
\end{displaymath}
and by \lemmaname~\ref{lemma-brenier} we have
\begin{displaymath}
\eta\left( \int_\setR h\left( u, v, w \right)\left( \xi \right) \d\xi \right) - \eta\left( 0 \right)
\le \int_\setR \eta^\prime\left( \xi \right) h\left( u, v, w \right)\left( \xi \right) \d\xi.
\end{displaymath}

Then using (\ref{definition-h}) and (\ref{definition-G}) on the right hand side, the cell entropy inequality (\ref{cell-entropy-inequality}) follows.
\end{proof}

\section{Main Result and Convergence of the Numerical Solutions}\label{main}

In this section we will present the main result of this paper, the convergence of the numerical solutions of the well-balanced scheme (see \definitionname~\ref{def-scheme}) to the entropy solution, and its proof.

\begin{theorem} \label{maintheorem}
Let (\ref{ass}), (\ref{assumptions-D}) and (\ref{engquist-osher flux}) be satisfied, let $( u_{\Delta x} )_{\Delta x > 0}$ be the sequence of numerical solutions of (\ref{ibvp}), in the sense of \definitionname~\ref{def-scheme} and assume the CFL-condition $\Lip_{K^{\Delta x}_T}\left( f \right) \frac{\Delta t}{\Delta x} \le 1$, where $C^{\Delta x}_T$ is defined as in (\ref{L-infty-bound}). Then there exists a $u \in L^\infty\left( \Omega_T \right)$, such that
\begin{equation} \intlabel{Aussage}
u_{\Delta x} \enspace \stackrel{\Delta x \to 0}{\longrightarrow} \enspace u
\quad \mbox{in} \enspace L^p\left( \Omega_T \right) \quad \mbox{for all } 1 \le p < \infty \quad
\end{equation}
and $u$ is the unique entropy solution of (\ref{ibvp}).
\end{theorem}

\begin{proof}
see Proof \ref{proof-maintheorem}
\end{proof}

For the proof we need the following lemmata.

\begin{lemma} \label{boundary-G}
Let $\left( \eta, q \right)$ be a boundary entropy pair, let $w \in \setR$ such that $\eta\left( w \right) = \eta^\prime\left( w \right) = q\left( w \right) = 0$ and let $G$ be defined by (\ref{definition-G}). Then we have for all $u, v \in \setR$
\begin{equation} \label{representation-G}
G\left( u, v \right)
= \int_w^u \eta^\prime\left( \xi \right) {f^\prime}^+\left( \xi \right) \d\xi
    - \int_w^v \eta^\prime\left( \xi \right) {f^\prime}^-\left( \xi \right) \d\xi.
\end{equation}

If additionally $u, v, w \in \left[ -C, C \right]$, $C > 0$, holds, the following inequality holds true:
\begin{equation} \label{boundary-inequality}
- \Lip_{\left[ -C, C \right]}\left( f \right) \eta\left( v \right)
\le G\left( u, v \right)
\le \Lip_{\left[ -C, C \right]}\left( f \right) \eta\left( u \right).
\end{equation}
\end{lemma}

\begin{proof}
Since $q\left( w \right) = 0$ we have
\begin{eqnarray*}
q\left( 0 \right)
&=& \int_w^0 q^\prime\left( \xi \right)  \d\xi
= \int_w^0 \eta^\prime\left( \xi \right) f^\prime\left( \xi \right) \d\xi\\
&=& \int_w^0 \eta^\prime\left( \xi \right) {f^\prime}^+\left( \xi \right) \d\xi
  - \int_w^0 \eta^\prime\left( \xi \right) {f^\prime}^-\left( \xi \right) \d\xi.
\end{eqnarray*}
Taking (\ref{definition-G}) into account, we obtain (\ref{representation-G}).

Since $\eta \in C^2\left( \setR \right)$ is convex and $\eta^\prime\left( w \right) = 0$, we have
\begin{displaymath}
\sgn\left( u - w \right) \eta^\prime\left( \xi \right) \ge 0 \quad \mbox{for all} \enspace \xi \in \left[ u, w \right].
\end{displaymath}
Therefore, as ${f^\prime}^+ \ge 0$, we see that
\begin{displaymath}
\int_w^u \eta^\prime\left( \xi \right) {f^\prime}^+\left( \xi \right) \d\xi
= \int_{[w, u]} \sgn\left( u - w \right) \eta^\prime\left( \xi \right) {f^\prime}^+\left( \xi \right) \d\xi
\ge 0.
\end{displaymath}
On the other hand $\eta\left( w \right) = 0$ and we can estimate:
\begin{eqnarray*}
\int_w^u \eta^\prime\left( \xi \right) {f^\prime}^+\left( \xi \right) \d\xi
&\le& \left\Vert f^\prime \right\Vert_{L^\infty\left( \left[ -C, C \right] \right)} \int_{[w, u]} \sgn\left( u - w \right) \eta^\prime\left( \xi \right) \d\xi\\
&=& \Lip_{\left[ -C, C \right]}\left( f \right) \eta\left( u \right).
\end{eqnarray*}

Applying the same arguments we see that
\begin{displaymath}
0 \le \int_w^v \eta^\prime\left( \xi \right) {f^\prime}^-\left( \xi \right) \d\xi
\le \Lip_{\left[ -C, C \right]}\left( f \right) \eta\left( v \right)
\end{displaymath}
and using (\ref{representation-G}) we obtain (\ref{boundary-inequality}).
\end{proof}

\begin{lemma}\label{lemma-discrete-process-sol}
Let (\ref{ass}), (\ref{assumptions-D}), (\ref{engquist-osher flux}) and the cell entropy inequality (\ref{cell-entropy-inequality}) be satisfied. For $\Delta t, \Delta x > 0$ assume that $\frac{\Delta t}{\Delta x} \le \lambda$ and let $( u^n_j )_{0 \le n \le N_T, j \in J}$ be the numerical solution of (\ref{ibvp}) in the sense of \definitionname~\ref{def-scheme}. Let $\left( \eta, q \right)$ be a boundary entropy pair, $w \in \setR$ such that $\eta\left( w \right) = \eta^\prime\left( w \right) = q\left( w \right) = 0$, $\vert u^n_j \vert, \vert w \vert \le C$ for all $j \in J$, $n \le N_T$ and $\vert u^n_{j+1,-} \vert, \vert u^n_{j-1,+} \vert \le C$ for all $j \in J$, $n < N_T$. Then for all $\varphi \in C_0^2( \setR^2 ) \cap C_0^\infty( \overline\Omega \times [ 0, T [ )$, $\varphi \ge 0$, we have
\begin{equation} \label{discrete-process-sol} \intlabel{result}
\begin{array}{c}
\displaystyle
\int_0^T \int_\Omega \eta\left( u_{\Delta x} \right) \partial_t \varphi + q\left( u_{\Delta x} \right) \partial_x \varphi
    - \eta^\prime\left( u_{\Delta x} \right) b\left( u_{\Delta x} \right) z^\prime\left( x \right) \varphi \d x \d t\next
\displaystyle
    + \int_\Omega \eta\left( u_0\left( x \right) \right) \varphi\left( x, 0 \right) \d x
    + \Lip_{\left[ -C, C \right]}\left( f \right) \int_0^T \eta\left( u_r\left( t \right) \right) \varphi\left( x_{j_r}, t \right) \d t\next
\displaystyle
    + \Lip_{\left[ -C, C \right]}\left( f \right) \int_0^T \eta\left( u_l\left( t \right) \right) \varphi\left( x_{j_l}, t \right) \d t
\ge \displaystyle \mathcal{O}\left( \Delta x \right).
\end{array}
\end{equation}
\end{lemma}

\begin{proof}
Multiplying the cell entropy inequality (\ref{cell-entropy-inequality}) by $\Delta x \, \varphi^n_j$ and summing over $n < N_T$ and $j \in J$ we obtain
\begin{equation} \intlabel{mainterm}
\begin{array}{rcl}
\displaystyle 0 &\ge& \displaystyle
\Delta x \sum_{n=0}^{N_T-1} \sum_{j \in J} \left( \eta\left( u^{n+1}_j \right) - \eta\left( u^n_j \right) \right) \varphi^n_j\next
&&\displaystyle
+ \Delta t \sum_{n=0}^{N_T-1} \sum_{j \in J} \left( G\left( u^n_j, u^n_{j+1} \right) - G\left( u^n_{j-1}, u^n_j \right) \right) \varphi^n_j\next
&&\displaystyle
+ \Delta t \sum_{n=0}^{N_T-1} \sum_{j \in J} \left( G\left( u^n_j, u^n_{j+1,-} \right) - G\left( u^n_j, u^n_{j+1} \right) \right) \varphi^n_j\next
&&\displaystyle
- \Delta t \sum_{n=0}^{N_T-1} \sum_{j \in J} \left( G\left( u^n_{j-1,+}, u^n_j \right) - G\left( u^n_{j-1}, u^n_j \right) \right) \varphi^n_j.
\end{array}
\end{equation}

We will now treat the four sums in (\intref{mainterm}) seperately. The first term is handled in a standard way by an index shifting argument:
\begin{equation} \intlabel{pre-term1}
\begin{array}{rcl}
\displaystyle \lefteqn{\sum_{n=0}^{N_T-1} \left( \eta\left( u^{n+1}_j \right) - \eta\left( u^n_j \right) \right) \varphi^n_j}\next
&=&\displaystyle
  - \sum_{n=0}^{N_T-1} \eta\left( u^n_j \right) \, \left( \varphi^n_j - \varphi^{n-1}_j \right)
  + \eta\left( u^{N_T}_j \right) \varphi^{N_T-1}_j
  - \eta\left( u^0_j \right) \varphi^{-1}_j\next
&=& \displaystyle
  - \Delta t \sum_{n=0}^{N_T-1} \eta\left( u^n_j \right)\, \left( \partial_t \varphi \right)^n_j
  - \eta\left( u^0_j \right) \varphi^0_j
  + \mathcal{E}_{t, j},
\end{array}
\end{equation}
where 
\begin{math}
\left\vert \mathcal{E}_{t, j} \right\vert
\le \Delta t \left\Vert \eta \right\Vert_{L^\infty\left( \left[ -C, C \right] \right)}
    \left( T \left\Vert \partial_t^2 \varphi \right\Vert_\infty + 2 \left\Vert \partial_t \varphi \right\Vert_\infty \right).
\end{math}

Using Jensen's inequality we deduce that
\begin{displaymath}
\eta\left( u^0_j \right)
\le \frac{1}{\Delta x} \int_{C_j} \eta\left( u_0\left( x \right) \right) \d x
\end{displaymath}
and from (\intref{pre-term1}) we obtain
\begin{equation} \intlabel{term1}
\begin{array}{l}
\displaystyle \lefteqn{\sum_{n=0}^{N_T-1} \left( \eta\left( u^{n+1}_j \right) - \eta\left( u^n_j \right) \right) \varphi^n_j}\next
\ge \displaystyle - \Delta t \sum_{n=0}^{N_T-1} \eta\left( u^n_j \right) \left( \partial_t \varphi \right)^n_j
    - \frac{1}{\Delta x} \int_{C_j} \eta\left( u_0\left( x \right) \right) \d x \varphi^0_j
    + \mathcal{E}_{t, j}.
\end{array}
\end{equation}

Now we turn to the second term of (\intref{mainterm}). Using (\ref{representation-G}), an index shifting argument and (\ref{boundary-inequality}) we obtain
\begin{eqnarray*}
\lefteqn{\sum_{j \in J} \left( G\left( u^n_j, u^n_{j+1} \right) - G\left( u^n_{j-1}, u^n_j \right) \right) \varphi^n_j}\\
&=& - \Delta x \sum_{j \in J} \int_w^{u^n_j} \eta^\prime\left( \xi \right) f^\prime\left( \xi \right) \d\xi \left( \partial_x \varphi \right)^n_j\\
&&  + G\left( u^n_{j_r}, u^n_{j_r+1} \right) \varphi^n_{j_r}
    - G\left( u^n_{j_l-1}, u^n_{j_l} \right) \varphi^n_{j_l}
    + \mathcal{E}_{x, n}\\
&\ge& - \Delta x \sum_{j \in J} q\left( u^n_j \right) \left( \partial_x \varphi \right)^n_j
    - \eta\left( u^n_{j_r+1} \right) \Lip_{\left[ -C, C \right]}\left( f \right) \varphi^n_{j_r}\\
&&  - \eta\left( u^n_{j_l-1} \right) \Lip_{\left[ -C, C \right]}\left( f \right) \varphi^n_{j_l}
    + \mathcal{E}_{x, n}.
\end{eqnarray*}
where $\vert \mathcal{E}_{x,n} \vert \le c_{x,n} \Delta x$ with
\begin{displaymath}
c_{x,n} = 2 C \left\Vert \eta^\prime \right\Vert_{L^\infty\left( \left[ -C, C \right] \right)}
    \left\Vert f^\prime \right\Vert_{L^\infty\left( \left[ -C, C \right] \right)}
    \left( \left\vert \Omega \right\vert \left\Vert \partial_x^2 \varphi \right\Vert_\infty
    + \left\Vert \partial_x \varphi \right\Vert_\infty \right).
\end{displaymath}

Again using Jensen's inequality to get
\begin{displaymath}
\eta\left( u^n_l \right) \le {1 \over \Delta t} \int_{t^n}^{t^{n+1}} \eta\left( u_l\left( t \right) \right) \d t,
\quad \eta\left( u^n_r \right) \le {1 \over \Delta t} \int_{t^n}^{t^{n+1}} \eta\left( u_r\left( t \right) \right) \d t,
\end{displaymath}
we arrive at the following inequality:
\begin{equation} \intlabel{term2}
\begin{array}{rcl}
\displaystyle
\lefteqn{\sum_{j \in J} \left( G\left( u^n_j, u^n_{j+1} \right) - G\left( u^n_{j-1}, u^n_j \right) \right) \varphi^n_j}\next
&\ge& \displaystyle
- \Delta x \sum_{j \in J} q\left( u^n_j \right) \left( \partial_x \varphi \right)^n_j
    - \Lip_{\left[ -C, C \right]}\left( f \right) \frac{1}{\Delta t} \int_{t^n}^{t^{n+1}} \eta\left( u_r\left( t \right) \right) \varphi^n_{j_r} \d t\next
&&\displaystyle
    - \Lip_{\left[ -C, C \right]}\left( f \right) \frac{1}{\Delta t} \int_{t^n}^{t^{n+1}} \eta\left( u_l\left( t \right) \right) \varphi^n_{j_l} \d t
    + \mathcal{E}_{x, n}.
\end{array}
\end{equation}

Now we turn our attention to the discretisation of the source term. By the mean value theorem there is a $\zeta^n_{j+1} \in [ u^n_{j+1}, u^n_{j+1,-} ]$ such that
\begin{eqnarray*}
\lefteqn{G\left( u^n_j, u^n_{j+1,-} \right) - G\left( u^n_j, u^n_{j+1} \right)}\\
&=& - \int_{u^n_{j+1}}^{u^n_{j+1,-}} \eta^\prime\left( \xi \right) {f^\prime}^-\left( \xi \right) \d\xi\\
&=& - \eta^\prime\left( \zeta^n_{j+1} \right) {f^\prime}^-\left( \zeta^n_{j+1} \right) \left( u^n_{j+1,-} - u^n_{j+1} \right).
\end{eqnarray*}
Also by the mean value theorem and (\ref{discrete equilibrium-states}) there is a $\vartheta^n_{j+1} \in [ u^n_{j+1}, u^n_{j+1,-} ]$ so that
\begin{displaymath}
z_{j+1} - z_j
= D\left( u^n_{j+1,-} \right) - D\left( u^n_{j+1} \right)
= D^\prime\left( \vartheta^n_{j+1} \right) \left( u^n_{j+1,-} - u^n_{j+1} \right)
\end{displaymath}
and since $\inf_\setR D^\prime > 0$ we have
\begin{equation} \intlabel{source1}
G\left( u^n_j, u^n_{j+1,-} \right) - G\left( u^n_j, u^n_{j+1} \right)
= - \Delta x \eta^\prime\left( \zeta^n_{j+1} \right) \frac{{f^\prime}^-( \zeta^n_{j+1} )}{D^\prime( \vartheta^n_{j+1} )}
    \frac{z_{j+1} - z_j}{\Delta x}.
\end{equation}

Since $z_{j_l} - z_{j_l-1} = z_{j_r+1} - z_{j_r} = 0$, multiplying (\intref{source1}) by $\varphi^n_j$ and summing over $j \in J$ we obtain
\begin{equation} \intlabel{term3}
\begin{array}{rcl}
\displaystyle \lefteqn{\sum_{j \in J} \left( G\left( u^n_j, u^n_{j+1,-} \right)  - G\left( u^n_j, u^n_{j+1} \right) \right) \varphi^n_j}\next
&=& \displaystyle - \Delta x \sum_{j \in J} \eta^\prime\left( \zeta^n_j \right)
    \frac{{f^\prime}^-( \zeta^n_j )}{D^\prime( \vartheta^n_j )} \frac{z_j - z_{j-1}}{\Delta x} \varphi^n_{j-1}\next
&=& \displaystyle - \Delta x \sum_{j \in J} \eta^\prime\left( u^n_j \right) \frac{{f^\prime}^-( u^n_j )}{D^\prime( u^n_j )}
    \frac{z_j - z_{j-1}}{\Delta x} \varphi^n_j + \mathcal{E}_{+,n},
\end{array}
\end{equation}
where
\begin{displaymath}
\mathcal{E}_{+,n} := \sum_{j \in J} \left[ \eta^\prime\left( u^n_j \right) \frac{{f^\prime}^-( u^n_j )}{D^\prime( u^n_j )}
    \varphi^n_j - \eta^\prime\left( \zeta^n_j \right) \frac{{f^\prime}^-( \zeta^n_j )}{D^\prime( \vartheta^n_j )} \varphi^n_{j-1} \right]
    (z_j - z_{j-1}).
\end{displaymath}

By the assumptions (\ref{assumptions-D}) $\frac{1}{D^\prime}$ is differentiable and we can estimate:
\begin{displaymath}
\begin{array}{l}
\displaystyle \left\vert \eta^\prime\left( \zeta^n_j \right) \frac{{f^\prime}^-( \zeta^n_j )}{D^\prime( \vartheta^n_j )}
    - \eta^\prime\left( u^n_j \right) \frac{{f^\prime}^-( u^n_j )}{D^\prime( u^n_j )} \right\vert\next
\le \displaystyle C^\ast \max\left\lbrace \left\vert \zeta^n_j - u^n_j \right\vert, \left\vert \vartheta^n_j - u^n_j \right\vert \right\rbrace
\le C^\ast \Delta x \frac{\left\Vert z^\prime \right\Vert_\infty}{\inf_\setR D^\prime}
\end{array}
\end{displaymath}
where $C^\ast = C^\ast\left( C, \eta, f, D \right)$. Using this result we see that $\vert \mathcal{E}_{+,n} \vert \le c_{\pm,n} \Delta x$ with
\begin{eqnarray*}
c_{\pm,n}
&=& \left\vert \Omega \right\vert C^\ast \frac{\left\Vert z^\prime \right\Vert^2_\infty}{\inf_\setR D^\prime}
    \left\Vert \varphi \right\Vert_\infty\\
&& + \left\vert \Omega \right\vert \frac{\left\Vert \eta^\prime \right\Vert_{L^\infty\left( \left[ -C, C \right] \right)}
      \left\Vert f^\prime \right\Vert_{L^\infty\left( \left[ -C, C \right] \right)}}{\inf_\setR D^\prime}
      \left\Vert z^\prime \right\Vert_\infty \left\Vert \varphi^\prime \right\Vert_\infty.
\end{eqnarray*}

In the same way we treat the other part of the source term and obtain
\begin{equation} \intlabel{term4}
\begin{array}{rcl}
\displaystyle \lefteqn{- \sum_{j \in J} \left( G\left( u^n_{j-1,+}, u^n_j \right) - G\left( u^n_{j-1}, u^n_j \right) \right) \varphi^n_j}\next
&=& \displaystyle \Delta x \sum_{j \in J} \eta^\prime\left( u^n_j \right) \frac{{f^\prime}^+\left( u^n_j \right)}{D^\prime\left( u^n_j \right)}
    {z_{j+1} - z_j \over \Delta x} \varphi^n_j + \mathcal{E}_{-,n},
\end{array}
\end{equation}
where $\vert \mathcal{E}_{-,n} \vert \le c_{\pm,n} \Delta x$.

Using (\intref{term1}), (\intref{term2}), (\intref{term3}) and (\intref{term4}) in (\intref{mainterm}) we arrive at the following inequality:
\begin{equation} \intlabel{pre-result}
\begin{array}{c}
\displaystyle
\Delta t \Delta x \sum_{n=0}^{N_T-1} \sum_{j \in J} \eta\left( u^n_j \right) \left( \partial_t \varphi \right)^n_j
+ \Delta t \Delta x \sum_{n=0}^{N_T-1} \sum_{j \in J} q\left( u^n_j \right) \left( \partial_x \varphi \right)^n_j\next
\displaystyle
- \Delta t \Delta x \sum_{0 \le n < N_T \above0pt j \in J} \eta^\prime\left( u^n_j \right)
    \left[  \frac{{f^\prime}^+( u^n_j )}{D^\prime( u^n_j )} \frac{z_{j+1} - z_j}{\Delta x}
      - \frac{{f^\prime}^-( u^n_j )}{D^\prime( u^n_j )} \frac{z_j - z_{j-1}}{\Delta x}
    \right] \varphi^n_j\next
\displaystyle
  + \sum_{j \in J} \int_{C_j} \eta\left( u_0\left( x \right) \right) \varphi^0_j \d x\next
\displaystyle
  + \Lip_{\left[ -C, C \right]}\left( f \right) \sum_{n=0}^{N_T-1}
        \int_{t^n}^{t^{n+1}} \eta\left( u_r\left( t \right) \right) \varphi^n_{j_r}
      +  \eta\left( u_l\left( t \right) \right) \varphi^n_{j_l} \d t
\ge \mathcal{O}\left( \Delta x \right).
\end{array}
\end{equation}

It is not hard to see that (\intref{pre-result}) implies (\intref{result}).
\end{proof}

For the final convergence proof we need the following relation between weak$\ast$ and strong convergence:

\begin{lemma} \label{L-p-convergence}
Let $\Omega$ be a bounded, measurable subset of $\setR^n$, let $\left( v_n \right)_{n \in \setN}$ be a bounded sequence in $L^\infty\left( \Omega \right)$ and $v \in L^\infty\left( \Omega \right)$ such that for all continuous functions $g \in C\left( \setR \right)$ we have
\begin{displaymath}
g\left( v_n \right) \to g\left( v \right) \quad \mbox{weak}\ast \enspace \mbox{in} \enspace L^\infty\left( \Omega \right).
\end{displaymath}
Then the convergence $v_n \to v$ is strong in $L^p\left( \Omega \right)$ for $1 \le p < \infty$.
\end{lemma}

\begin{proof}
see \cite{vovelle}, Lemma~10.
\end{proof}

\begin{numproof}[Proof of \theoremname~\ref{maintheorem}] \label{proof-maintheorem}
Let $\delta > 0$ be fixed. Then for $\Delta x \le \delta$ we have $C^{\Delta x}_T \le C^\delta_T$, where $C^{\Delta x}_T$ is defined by (\ref{L-infty-bound}). Hence the CFL-condition and \lemmaname~\ref{lemma-boundedness} provide for $\Delta x \le \delta$
\begin{displaymath}
\left\Vert u_{\Delta x} \right\Vert_{L^\infty\left( \Omega \times \left] 0, T \right[ \right)}
\le C^\delta_T.
\end{displaymath}
This lemma also provides the boundedness of the discrete equilibrium states, as defined in (\ref{discrete equilibrium-states}).

Let $\left( \eta, q \right)$ be a boundary entropy pair and let $w \in \setR$ such that
\begin{math}
\eta\left( w \right) = \eta^\prime\left( w \right) = q\left( w \right) = 0.
\end{math}
By \lemmaname~\ref{lemma-cell-entropy-inequality} there is a consistent numerical entropy flux $G$, so that the cell entropy inequality (\ref{cell-entropy-inequality}) is satisfied.

Let $\varphi \in C_0^\infty( \overline\Omega \times [ 0, T [ )$, $\varphi \ge 0$. Then $\varphi$ can be extended to a function in $C_0^\infty( \overline\Omega \times [ 0, T [ ) \cap C_0^2( \setR^2 )$, which we will also denote as $\varphi$.

Setting $C := \max\lbrace C^\delta_T, \vert w \vert \rbrace$ all the premises of \lemmaname~\ref{lemma-discrete-process-sol} are fulfilled and we can deduce that for $0 < \Delta x \le \delta$ and for all non-negative test functions $\varphi \in C_0^\infty( \overline \Omega \times [ 0, T [ ) \cap C_0^2( \setR^2 )$, we have (\ref{discrete-process-sol}).

Now \lemmaname~\ref{lemma-eymard} allows us to find a sequence $( \Delta x_k )_{k \in \setN} \subset ] 0, \delta ]$ and a
$\mu \in L^\infty( \Omega_T \times ] 0, 1 [ )$ such that for all $g \in C( \setR )$ we have
\begin{equation} \intlabel{weak-convergence}
g\left( u_{\Delta x_k} \right) \to \int_0^1 g\left( \mu\left( \cdot, \cdot, \alpha \right) \right) \d\alpha
\quad \mbox{weak}\ast \enspace \mbox{in}\enspace L^\infty\left( \Omega_T \right).
\end{equation}

Since $\eta, q, \eta^\prime b \in C( \setR )$ and  $\partial_t \varphi, \partial_x \varphi, z^\prime \varphi \in L^1( \Omega_T )$ we can pass to the limit in (\ref{discrete-process-sol}) and obtain (\ref{process-sol}). Thus $\mu$ is an entropy process solution of (\ref{ibvp}) and by \theoremname~\ref{uniqueness} get a $u \in L^\infty( \Omega_T )$ such that
\begin{math}
\mu( \cdot, \cdot, \alpha ) = u
\end{math}
a.e.\ for almost all $\alpha \in ] 0, 1 [$.

Since the entropy process solution is unique by \theoremname~\ref{uniqueness}, we conclude that any subsequence of $( u_{\Delta x} )_{\Delta x > 0}$ converges to weak$\ast$ to $u$ as $\Delta x \to 0$. By \lemmaname~\ref{L-p-convergence} we see that $( u_{\Delta x} )_{\Delta x > 0}$ converges strongly to $u$ in $L^p( \Omega_T )$.
\end{numproof}

\section{Numerical Experiments}\label{numerics}

In this section we will numerically test the well-balanced scheme (see \definitionname~\ref{def-scheme}) and compare it to the standard discretization (\ref{standard-scheme}). Therefore we take a look at the Burgers-Hopf equation
\begin{equation} \label{burgers-hopf}
\partial_t u + \partial_x \frac{u^2}{2} + z^\prime\left( x \right) u = 0
\end{equation}
on $\Omega := ] 0, 4 [$, which means $D( s ) = s$. So (\ref{assumptions-D}) is obviously satisfied. For the source term we choose
\begin{equation} \label{example-z}
z\left( x \right) = \left\lbrace
\begin{array}{ll}
\cos\left( \pi x \right) & \enspace \mbox{for} \enspace x \in \left( \frac{3}{2}, \frac{5}{2} \right),\\
0 & \enspace otherwise
\end{array}\right..
\end{equation}
Notice that $z^\prime$ is discontinuous and $u_c( x ) := c - z( x )$ is an equilibrium of (\ref{burgers-hopf}).

\begin{testcase} \label{test-equilibrium}
Assume the equilibrium initial data $u_2 := 2 - z( x )$ and the constant boundary data $u_l = u_r = 2$. The exact solution to this problem is given by $u( x, t ) = u_2( x )$.

Now we use 40 nodes in space and compare the numerical solutions of both schemes (as defined in (\ref{def-scheme}) and (\ref{standard-scheme}) respectively) at $t = 3$. The result is shown in \figurename~\ref{figure-equilibrium}.
\begin{figure}
\centering\resizebox{0.9\textwidth}{!}{\includegraphics{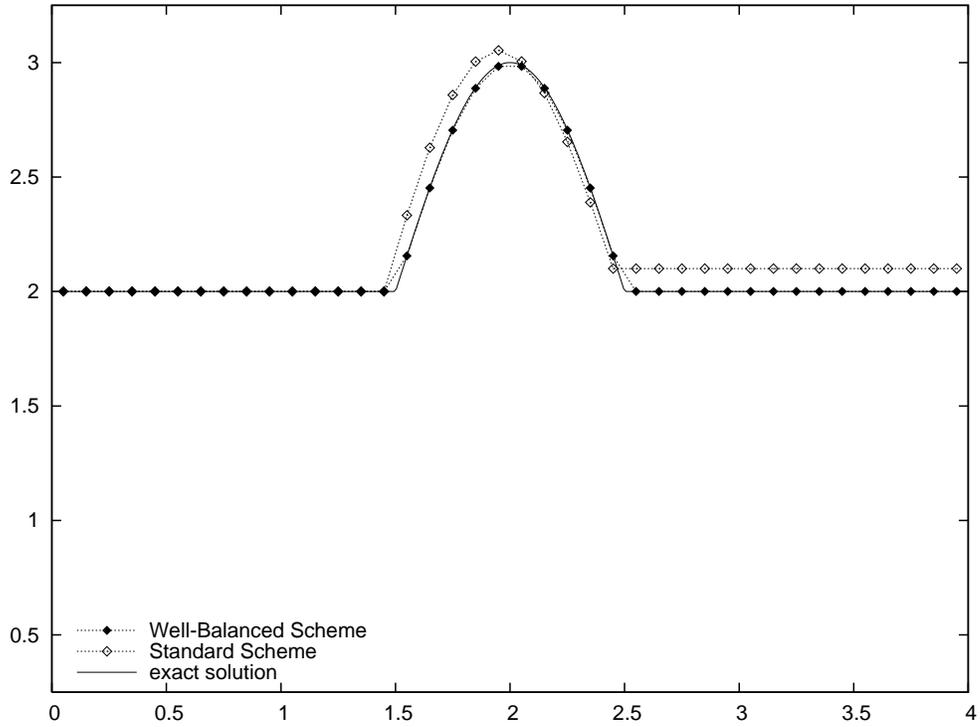}}
\caption{Well-Balanced and Standard Scheme for Test~Case~\ref{test-equilibrium} at $t=3$}
\label{figure-equilibrium}
\end{figure}

The \tablename s~\ref{table-equilibrium-balanced} and \ref{table-equilibrium-standard} show the $L^1$-error for the numerical solution obtained for different grid sizes $\Delta x$. They also display the $L^1$-error of the numerical solution with respect to the projection of the exact solution to the space of piecewise constant functions (Numerical Error).
\begin{table}
\begin{center}\begin{tabular}{rrrrr}
\multicolumn{1}{c}{$\Delta x$} & \multicolumn{1}{c}{$\Delta t$}
& \multicolumn{1}{c}{$L^1$-Error} & \multicolumn{1}{c}{Num.\ Error}
& \multicolumn{1}{c}{CPU-Time}\\
\hline\noalign{\smallskip}
$10^{-1}$ & $5.7 \cdot 10^{-6}$ & $5.02 \cdot 10^{-2}$ & $0$ & $0.62 \, s$\\
$10^{-2}$ & $3.5 \cdot 10^{-6}$ & $5.00 \cdot 10^{-3}$ & $0$ & $8.83 \, s$\\
$10^{-3}$ & $7.3 \cdot 10^{-7}$ & $5.00 \cdot 10^{-4}$ & $0$ & $423.2 \, s$\\
$10^{-4}$ & $8 \cdot 10^{-8}$ & $5.00 \cdot 10^{-5}$ & $0$ & $42\,688 \, s$\\
\noalign{\smallskip}\hline
\end{tabular}\end{center}
\caption{Well-Balanced Scheme for Equilibrium Initial Data $u_0 = z\left( x \right)$}
\label{table-equilibrium-balanced}
\end{table}

\begin{table}
\begin{center}\begin{tabular}{rrrrr}
\multicolumn{1}{c}{$\Delta x$} & \multicolumn{1}{c}{$\Delta t$}
& \multicolumn{1}{c}{$L^1$-Error} & \multicolumn{1}{c}{Num. Error}
& \multicolumn{1}{c}{CPU-Time}\\
\hline\noalign{\smallskip}
$2 \cdot 10^{-2}$ & $4.5 \cdot 10^{-6}$ & $5.07 \cdot 10^{-2}$ & $4.90 \cdot 10^{-2}$ & $2.94 \, s$\\
$2 \cdot 10^{-3}$ & $1.3 \cdot 10^{-6}$ & $5.09 \cdot 10^{-3}$ & $4.91 \cdot 10^{-3}$ & $97.6 \, s$\\
$2 \cdot 10^{-4}$ & $1.6 \cdot 10^{-7}$ & $5.09 \cdot 10^{-4}$ & $4.92 \cdot 10^{-4}$ & $10\,885 \, s$\\
$2 \cdot 10^{-5}$ & $2 \cdot 10^{-8}$ & $5.06 \cdot 10^{-5}$ & $4.97 \cdot 10^{-5}$ & \multicolumn{1}{c}{---} \\
\noalign{\smallskip}\hline
\end{tabular}\end{center}
\caption{Standard Scheme for Equilibrium Initial Data $2 - z\left( x \right)$}
\label{table-equilibrium-standard}
\end{table}
\end{testcase}

In these computations $\Delta t$ is computed such that the CFL-condition is satisfied. The necessary estimate on $\vert u^n_j \vert$ is given by \lemmaname~\ref{lemma-boundedness}. For the standard discretisation (\ref{standard-scheme}) we have
\begin{displaymath}
\left\vert u^{n+1}_j \right\vert
\le \sup \left\vert u^n_j \right\vert
    + \Delta t \left\Vert z^\prime \right\Vert_{L^\infty\left( \Omega \right)} \sup \left\vert u^n_j \right\vert
\le M e^{T \left\Vert z^\prime \right\Vert_{L^\infty\left( \Omega \right)}}
\le C^{\Delta x}_T,
\end{displaymath}
where $C^{\Delta x}_T$ is defined as in (\ref{L-infty-bound}). So our estimate for the well-balanced scheme seems far from optimal.

\begin{testcase} \label{test-cmp-equilibrium}
Assume the constant initial data $u_0 = 1$ and the constant boundary data $u_l = 2$, $u_r = 1$.Again we use $40$ nodes in space to compare the well-balanced scheme as defined in (\ref{def-scheme}) to the standard discretization (\ref{standard-scheme}). A reference solution was obtained by the standard scheme with $40\,000$ nodes in space. The result is displayed in \figurename~\ref{figure-cmp-equilibrium}.
\begin{figure}
\centering{\begin{minipage}{0.99\linewidth}
\begin{minipage}[t]{0.495\linewidth}
\centering\resizebox{\linewidth}{!}{\includegraphics{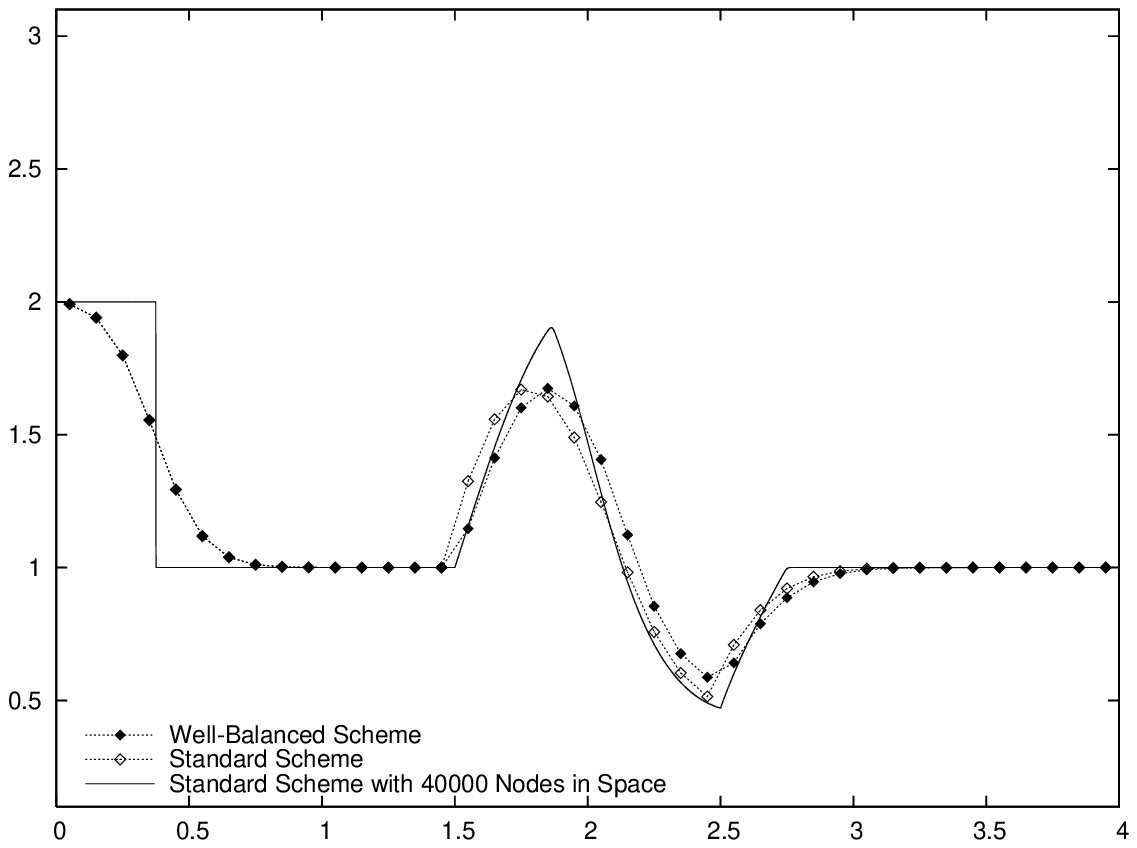}}\\$t = 0.25$
\end{minipage}
\hfill
\begin{minipage}[t]{0.495\linewidth}
\centering\resizebox{\linewidth}{!}{\includegraphics{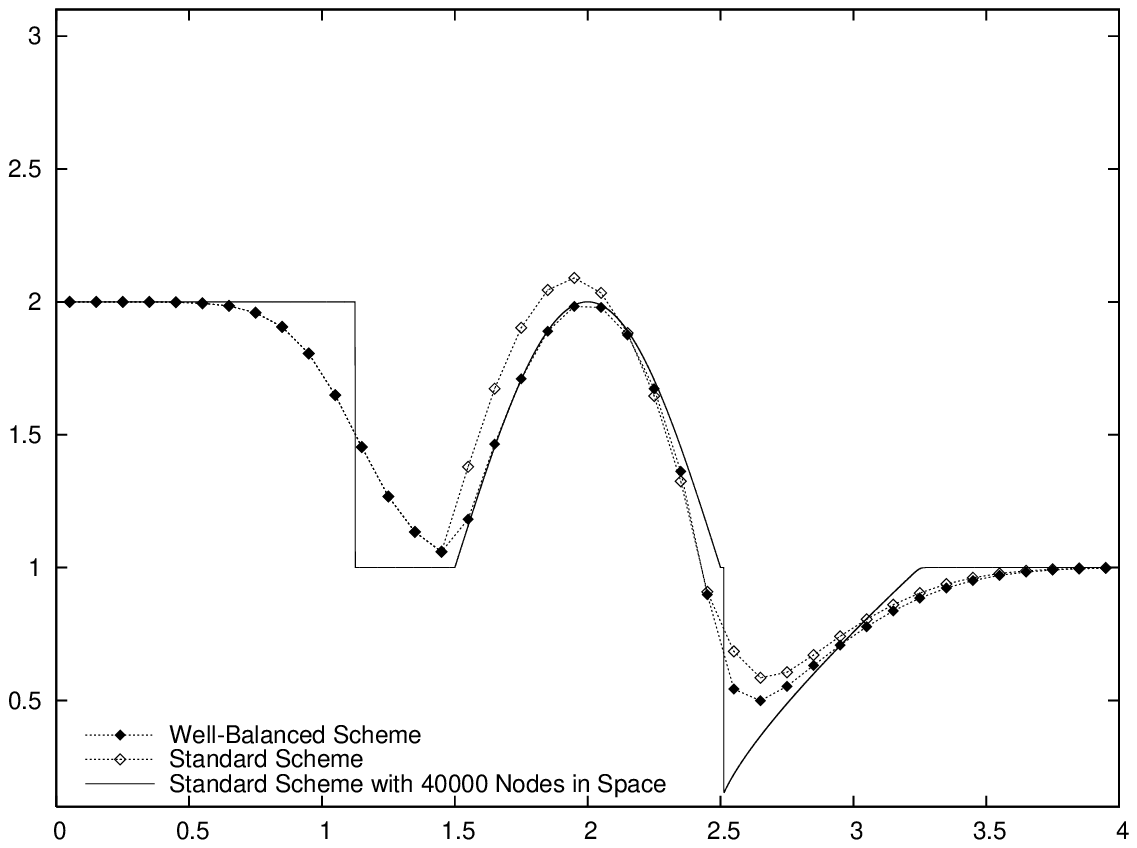}}\\$t = 0.75$
\end{minipage}
\par\bigskip\begin{minipage}[t]{0.495\linewidth}
\centering\resizebox{\linewidth}{!}{\includegraphics{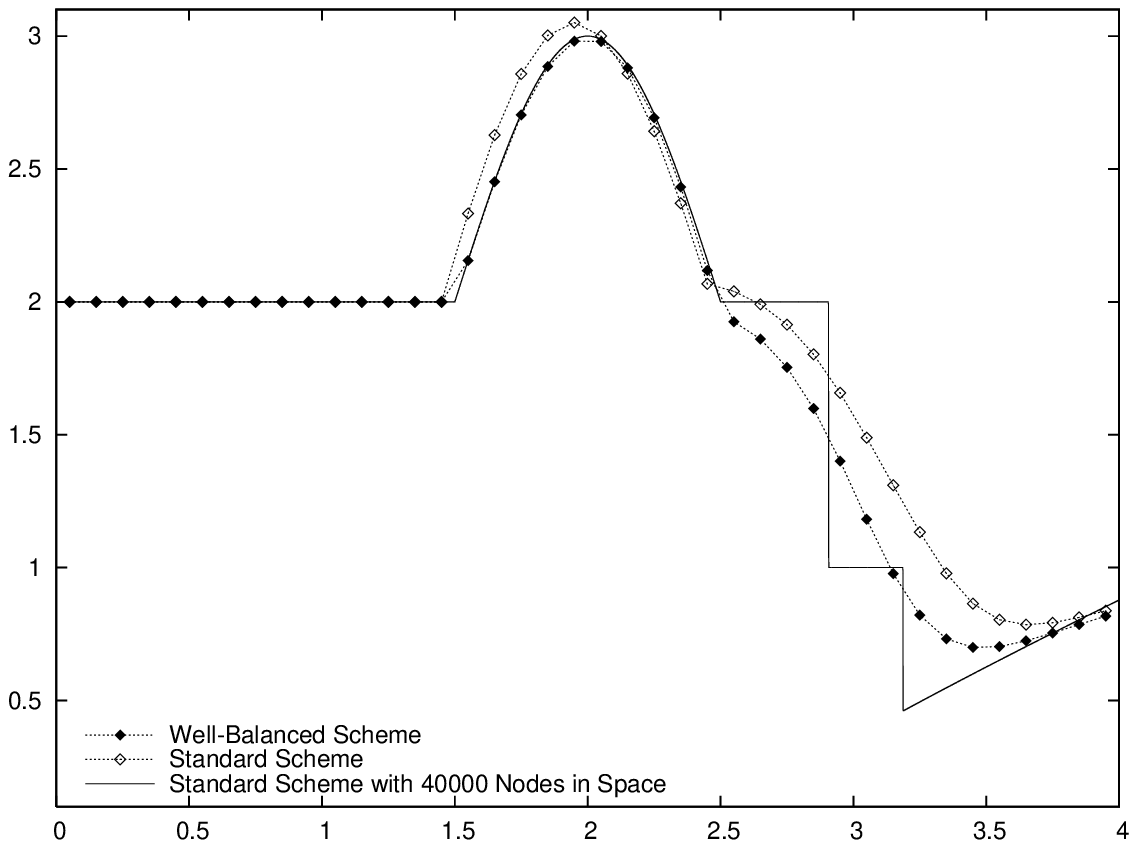}}\\$t = 1.75$
\end{minipage}
\hfill
\begin{minipage}[t]{0.495\linewidth}
\centering\resizebox{\linewidth}{!}{\includegraphics{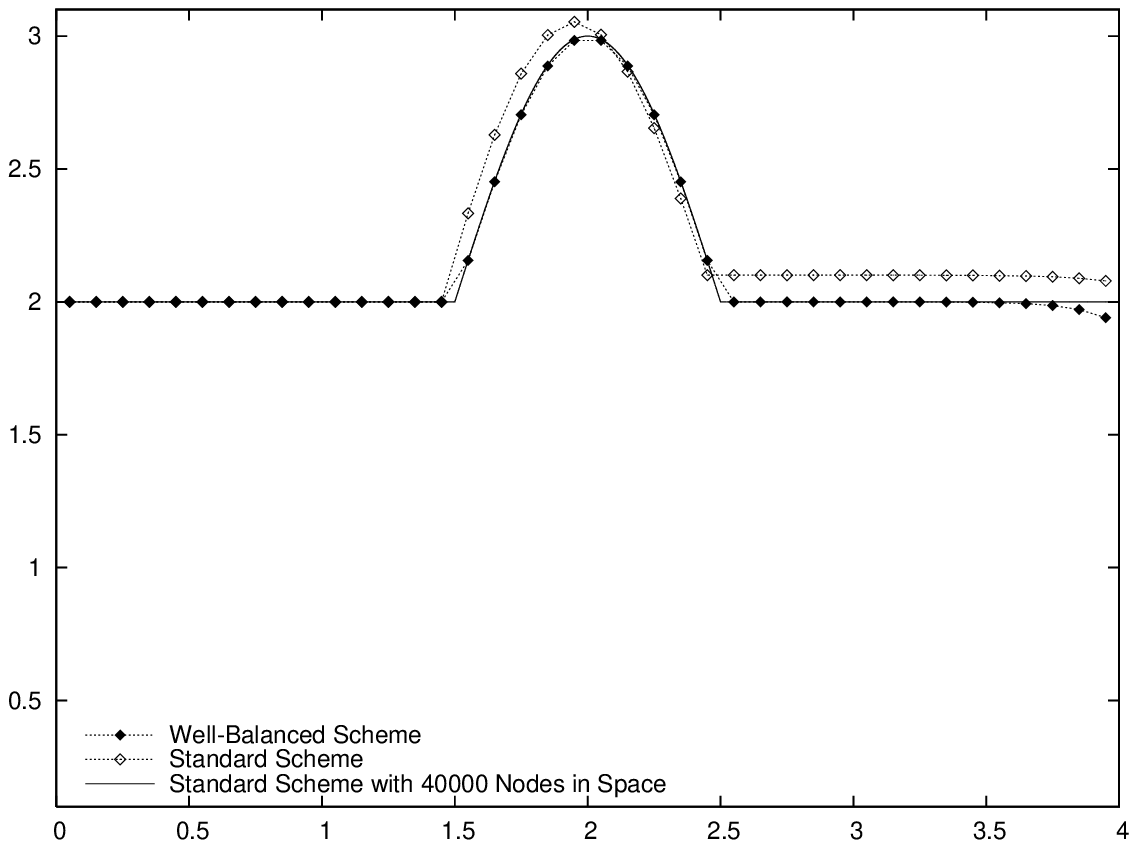}}\\$t = 2.75$
\end{minipage}
\caption{Well-balanced scheme vs.\ standard discretization for Test~Case~\ref{test-cmp-equilibrium}}
\end{minipage}}
\label{figure-cmp-equilibrium}
\end{figure}
\end{testcase}

Now we turn to the weakness of the well-balanced scheme. For the Burgers-Hopf equation (\ref{burgers-hopf}) $u( x, t ) = 0$ is an equilibrium for any choice of $z$. Thus, if we ask the boundary data $u_l = u_r = 0$ there are two equilibria satisfying this condition: $-z(x)$ and $0$.
 
\begin{testcase} \label{test-null}
We assume the constant initial data $u_0 = 0$ and the constant boundary data $u_l = u_r = 0$. the exact solution is given by $u\left( x, t \right) = 0$. \figurename~\ref{figure-null-equilibrium} displays the numerical solutions obtained by the well-balanced scheme for $40$, $400$, $4\,000$ and $40\,000$ nodes in space at $t=2.5$.
\begin{figure}
\centering\resizebox{0.9\textwidth}{!}{\includegraphics{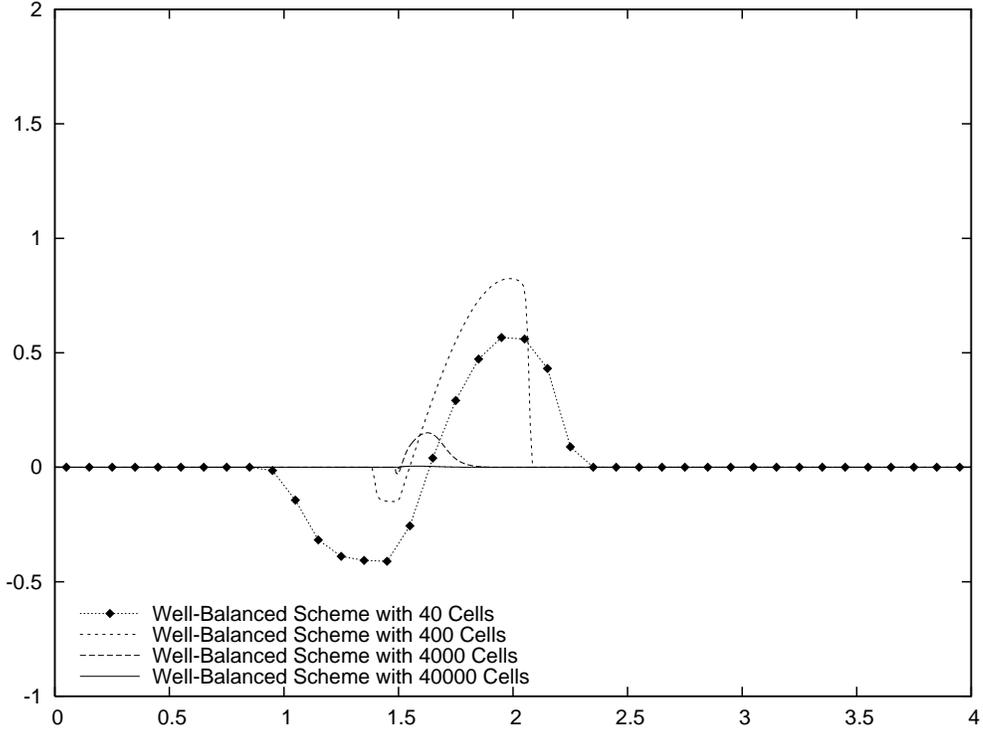}}
\caption{Well-balanced scheme for Test~Case~\ref{test-null} at $t=2.5$}
\label{figure-null-equilibrium}
\end{figure}

\tablename~\ref{table-null-equilibrium} shows the $L^1$-error for this problem, which is not acceptable even for $4\,000$ nodes in space.
\begin{table}
\begin{center}\begin{tabular}{rrrr}
\multicolumn{1}{c}{$\Delta x$} & \multicolumn{1}{c}{$\Delta t$}
& \multicolumn{1}{c}{$L^1$-Error} & \multicolumn{1}{c}{CPU-Time}\\
\hline\noalign{\smallskip}
$10^{-1}$ & $6.14 \cdot 10^{-6}$ & $4.388 \cdot 10^{-1} $ & $0.90 \, s$\\
$10^{-2}$ & $6.14 \cdot 10^{-6}$ & $3.164 \cdot 10^{-1} $ & $5.23 \, s$\\
$10^{-3}$ & $6.14 \cdot 10^{-6}$ & $2.678 \cdot 10^{-2} $ & $47.9 \, s$\\
$10^{-4}$ & $6.14 \cdot 10^{-6}$ & $8.421 \cdot 10^{-4} $ & $740.8 \, s$\\
\noalign{\smallskip}\hline
\end{tabular}\end{center}
\caption{Well-Balanced Scheme for Initial Data $u_0 = 0$}
\label{table-null-equilibrium}
\end{table}
\end{testcase}

Notice that the convergence proof for the well-balanced scheme holds for $u = 0$, even though we divided the differential equation in (\ref{ibvp}) by $u$ in the motivation.

We also tested both the well-balanced scheme and the standard scheme (\ref{standard-scheme}) for a discontinuous $z$:
\begin{testcase} \label{test-disc-equilibrium}
Let $z$ be given by the discontiuous function
\begin{displaymath}
z\left( x \right) = \left\lbrace
\begin{array}{ll}
\sin\left( \pi x \right) & \enspace \mbox{for} \enspace x \in \left( \frac{3}{2}, \frac{5}{2} \right),\\
0 & \enspace otherwise
\end{array}\right.
\end{displaymath}
and choose the initial data $u_2 = 2 - z( x )$ and the constant boundary data $u_l = u_r = 2$. The numerical results in this case are displayed in \figurename~\ref{figure-disc-equilibrium}.
\begin{figure}
\centering\resizebox{0.9\textwidth}{!}{\includegraphics{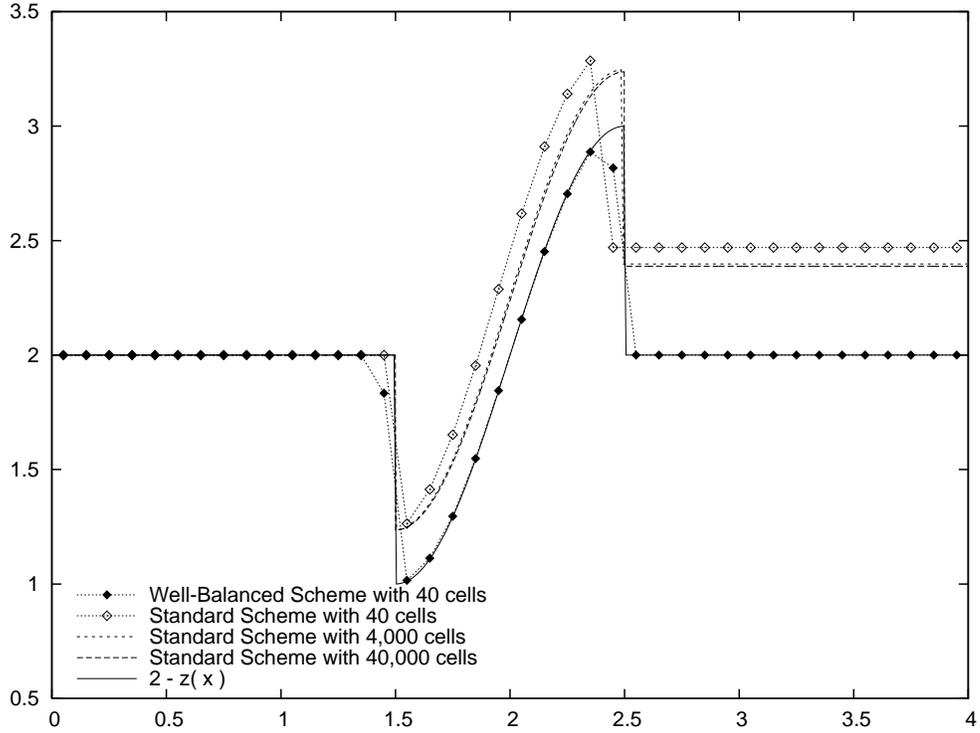}}
\caption{Well-Balanced and Standard Scheme for Test~Case~\ref{test-disc-equilibrium} at $t=3$}
\label{figure-disc-equilibrium}
\end{figure}

Note that while the well-balanced scheme converges to the equilibrium quite fast, the standard scheme seems not to converge to the equilibrium at all.
\end{testcase}

\end{document}